\documentclass[11pt]{article}
\usepackage{amsmath,amsthm,amssymb,amscd}
\usepackage{hyperref}
\usepackage{geometry}
\usepackage{tikz-cd}
\usepackage{booktabs} 
\usepackage{graphicx}
\usepackage{subcaption} 

\usepackage{multirow}
\usepackage{tikz}
\usetikzlibrary{calc,positioning,3d}
\usepackage[normalem]{ulem} 
\usepackage{imakeidx}
\usepackage{makecell}
\usepackage{float}
\makeindex
\usepackage{seqsplit}
\newcommand{\ttseq}[1]{{\ttfamily\seqsplit{#1}}}
\usepackage{xr-hyper}
\usepackage{cleveref}

\newtheorem{theorem}{Theorem}[section]

\theoremstyle{definition}

\newtheorem{example}[theorem]{Example}

\title{GBNL: Graded Betti Number Learning of Complex Biological Data}

\author{Mushal Zia$^1$,  
	Faisal Suwayyid $^{1,2}$, 
	and Guo-Wei Wei\footnote{Corresponding author: Guo-Wei Wei (weig@msu.edu).~}$^{~1,3,4}$ \\		
	$^1$Department of Mathematics,\\
	Michigan State University, MI 48824, USA.\\
	$^2$Department of Mathematics,\\
	King Fahd University of Petroleum and Minerals, Dhahran 31261, KSA.\\
	$^3$Department of Mathematics and Statistics,\\
	Michigan State University, MI 48824, USA.\\
	$^4$Department of Biochemistry and Molecular Biology,\\
	Michigan State University, MI 48824, USA.
}

\date{} 
\begin{document}
	\maketitle 
	
	\begin{abstract}	    
While persistent homology is widely used for data shape analysis, persistent commutative algebra (PCA) has seen limited adoption in machine learning and data science. Unlike persistent homology, which delivers topological invariants in the form of Betti numbers, PCA provides both algebraic invariants and graded Betti numbers. However, graded Betti numbers have seldom been applied to real-world data. In this work, we introduce the first-of-its-kind application of commutative algebra graded Betti numbers in machine learning and data science. Specifically, we present Graded Betti Number Learning (GBNL) for protein-nucleic acid binding prediction. Protein-DNA/RNA interactions are fundamental to cellular processes such as replication, transcription, translation, and gene regulation, and their understanding and prediction remain challenging. GBNL represents each nucleic acid sequence as a family of $k$-mer-specific sets and derives persistent graded Betti invariants from PCA, generating multiscale topological representations of local nucleotide organization. To incorporate cross-molecule context, these graded Betti representations are paired with transformer-based protein embeddings, linking nucleotide-level signals with global protein patterns. The proposed graded Betti representations effectively detect single-site mutations and distinguish complete mutation patterns. Operating on primary sequences with minimal preprocessing, GBNL bridges commutative algebra, reduced algebraic topology, combinatorics, and machine learning, establishing a new paradigm for comparative sequence analysis. Numerical studies using three datasets highlight the success of GBNL in protein-nucleic acid binding prediction.
		
	\end{abstract}
	
	
	\newpage	
	{\setcounter{tocdepth}{4} \tableofcontents}
	\setcounter{page}{1}
	\newpage	
	
	\section{Introduction}
	
	Protein-nucleic acid interactions underpin the molecular framework of gene regulation and flow of genetic information. Through sequence-specific recognition of DNA and RNA, proteins orchestrate replication, repair, transcription, splicing, and translation, thereby coupling genotype to controlled cellular function. Moreover, the interpretation and execution of genetic instructions inside a cell are correctly governed by the specificity and stability of these contacts. Quantification of protein-nucleic acid binding affinities is thus critical for not only understanding recognition principles but also for enabling synthetic biology, genome engineering, and therapeutic innovation~\cite{bertoldo2023rna}. At atomic resolution, affinity emerges from a calibrated balance of hydrogen bonds, electrostatics, hydrophobic packing, and shape complementarity. Disruptions in any of these processes can miswire regulation and are linked to neurodegeneration, autoimmunity, and cancer~\cite{kisby2024introducing}. Therefore, conserving both global organization and nucleotide-resolution details remains a persistent challenge molecular biophysics and computational biology.
	
	Standard binding assays such as electrophoretic mobility shift assays (EMSA), surface plasmon resonance (SPR), isothermal titration calorimetry (ITC), fluorescence readouts, and filter binding, deliver accurate affinities yet they remain costly and scale poorly~\cite{stockley2009filter}. Computational substitutes reduce cost however, they introduce tradeoffs. Moreover, physics-based free energy calculations such as thermodynamic integration, free energy perturbation, and molecular mechanics Poisson-Boltzmann surface area (MM-PBSA) yield detailed estimates but become prohibitive for large or flexible complexes. Furthermore, knowledge-based potentials \cite{zhang2005knowledge}, empirical and force-field scoring functions 
	\cite{nithin2019structure}, and machine learning with engineered descriptors \cite{bitencourt2019machine} improve speed but can lack transferability and interpretability. These approaches perform well for protein-protein and protein-ligand systems, yet often degrade on DNA due to nucleic acid conformational variability and limited curated affinities. For RNA, crosslinking immunoprecipitation sequencing \cite{hafner2010transcriptome}, sequence driven predictors \cite{zhao2011structure,liu2016prediction,shen2023svsbi}, knowledge-based scoring \cite{tuszynska2011dars}, and coarse grained docking \cite{setny2011coarse} broaden coverage, but quantitative affinity measurements remain sparse because of molecular flexibility and experimental constraints \cite{iwakiri2016improved}.
	
	Recent progress in computational molecular biology has been driven by data and machine learning models that now play a central role in predicting biomolecular interactions and accelerating drug discovery \cite{lin2023evolutionary,song2024multiobjective}. The integration of bioinformatics with modern learning frameworks has enabled large scale inference of complex physical and chemical relationships from high-dimensional molecular data \cite{lo2018machine}. Within this evolution, mathematical artificial intelligence (AI), particularly topological deep learning (TDL) introduced in 2017 \cite{cang2017topologynet}, has emerged as a rational approach for representing and learning biomolecular structure \cite{papamarkou2024position,nguyen2020review} and has demonstrated strong performance in community challenges such as the D3R Grand Challenge series \cite{nguyen2019mathematical,nguyen2020mathdl}. This framework has since inspired the development of topological sequence analysis for genomic modeling \cite{hozumi2024revealing}, where delta complex formulations allow efficient treatment of large sequence data \cite{liu2025topological2}, and category theoretic formulations improve resolution among closely related genetic variants \cite{liu2025topological}. In parallel, machine learning models for protein-DNA binding have evolved from early atomic pairwise statistical potentials \cite{zhao2010structure}, sequence-based affinity profiling \cite{rastogi2018accurate}, and physics-informed methods such as DNAffinity \cite{barissi2022dnaffinity} to ensemble architectures including PreDBA \cite{yang2020predba}, PDA-Pred \cite{harini2023pda}, and emPDBA \cite{yang2023empdba}. Following the first quantitative protein-RNA affinity dataset~\cite{yang2013dataset}, the trajectory moved toward structure-guided prediction with methods by Deng et al.~\cite{deng2019predprba} and Nithin et al.~\cite{nithin2019structure}, with subsequent refinements by Hong et al.~\cite{hong2023updated} and Harini et al.~\cite{harini2024pred}. However, despite these advances, most approaches continue to rely on stacked regressors, class-specific architectures, or interface-level descriptors to mitigate data scarcity, underscoring the need for frameworks that combine mathematical rigor with biological interpretability.
	
	Persistent homology (PH) serves as a foundational method in TDL for quantifying multiscale shape analysis via the birth and death of topological features. The resulting representations built from Betti numbers, barcodes, and Betti curves yield concise encodings of protein structure, molecular surfaces, and genomic organization~\cite{nguyen2020review,wee2025review}. Because it records only topological invariants, PH can miss local geometric or physicochemical changes. The persistent Laplacian (PL), a spectral approach introduced in 2019 \cite{wang2020persistent,chen2019evolutionary}, augments PH with non-harmonic spectra to incorporate geometric details and has progressed in theory and applications, including protein engineering and viral variant analysis \cite{memoli2022persistent,qiu2023persistent,chen2022persistent}. Nevertheless, both PH and spectral approaches primarily emphasize global information and connectivity.
	
	Commutative algebra offers a complementary route by translating simplicial and combinatorial data into polynomial ideals and algebraic invariants \cite{eisenbud2013commutative}. Persistent Stanley-Reisner theory integrates multiscale analysis with this algebraic setting to produce persistent graded Betti numbers, persistent facet ideals, $f$-vectors, $h$-vectors, as algebraic analogues of topological summaries \cite{suwayyid2025persistent}. These new representations follow the evolution of square-free monomial ideals under filtration, revealing how geometric and physicochemical interactions emerge or vanish across scales. The same algebraic machinery has already improved structure based protein-ligand affinity prediction through a commutative algebra representation of interfacial geometry and chemistry \cite{feng2025caml}. Extending this idea, a commutative algebra neural network learns graded algebraic signatures that trace mutation driven changes and relate them to disease mechanisms \cite{wee2025commutative}. In parallel, a $k$-mer formulation casts sequence positions into persistent algebraic features, enabling comparative genomics and phylogenetic inference at scales \cite{suwayyid2025cakl}. Together, these directions establish persistent commutative algebraic invariants as interpretable, data efficient descriptors that connect topology, geometry, and learning for biomolecular systems.
	
	In this study, we present GBNL, the first graded Betti framework integrating commutative algebra with biomolecular modeling and machine learning. It encodes nucleic acid organization through graded Betti numbers tracked over a positional filtration for protein-nucleic acid binding prediction. We pair these algebraic/topological descriptors with transformer-based protein embeddings to capture interface chemistry at two complementary scales. The result is a compact, interpretable feature set for predicting protein-DNA and protein-RNA binding. We validate the approach on three standard datasets and a targeted mutation analysis of the N-China-F primer of the SARS-CoV-2 N gene, showing that the graded Betti signal isolates local sequence changes that conventional Betti counts miss.

\section{Results and Discussion}\label{sec:results}

\subsection{Overview of the Workflow}

\begin{figure}[t!]
	\centering
	\begin{subfigure}[b]{1\textwidth}
		\includegraphics[width=\textwidth]{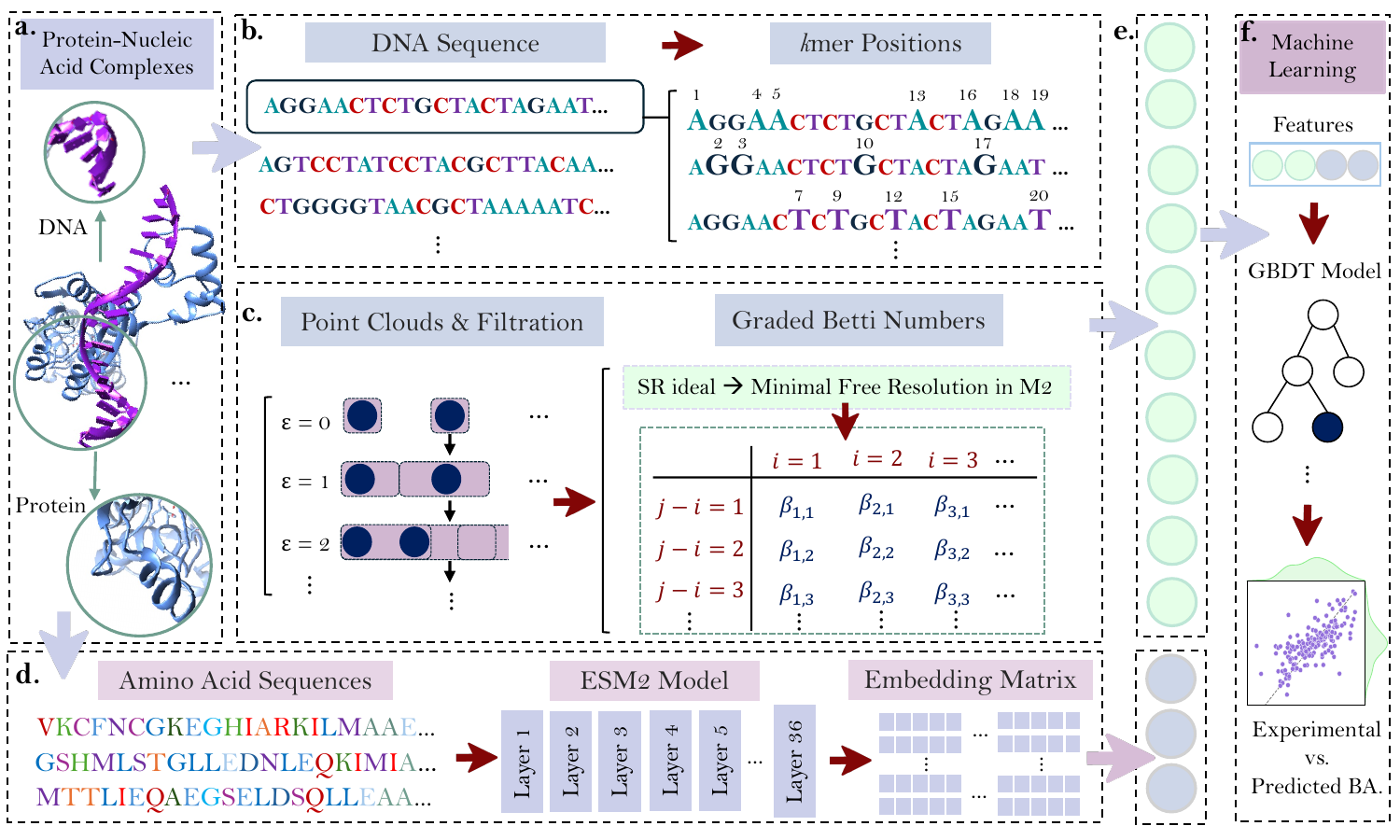}
	\end{subfigure}
	\caption{\textbf{An Illustration of GBNL workflow.}
		(a). Representative protein-nucleic acid complex provide biological context; downstream computations use only the DNA/RNA sequence (for per-nucleotide point clouds) and the protein amino-acid sequence (for embeddings).
		(b). \(k\)-mers are extracted from the DNA sequences of the given complexes. For each \(k=1\), the set of its occurrence positions within the sequence is treated as input data.
		(c). For each nucleotide, a Vietoris-Rips (VR) complex with maximum dimension 1 is constructed from its position-based point cloud across an \(\varepsilon\)-filtration. At each selected \(\varepsilon\), the maximal simplices define a simplicial complex whose Stanley-Reisner (SR) ideal \(I^{\varepsilon}\) is formed. The minimal free resolution of \(S/I^{\varepsilon}\), where \(S\) is a polynomial ring, is then computed in Macaulay2 (M2), yielding the graded Betti numbers \(\beta_{i,j}\), which serve as the topological features for the corresponding nucleotide. 
		(d). From the protein amino acid sequences of the same complexes, 2560-dimensional embedding vectors are generated using the state-of-the-art ESM2 model with 36 layers.
		(e). The extracted graded Betti features of DNA and the ESM2 embeddings of proteins are concatenated to form a unified algebraic-sequence representation.
		(f). Finally, a gradient boosted decision tree (GBDT) model is trained using these combined features for the prediction of experimental binding affinities.}
	\label{fig:workflow}
\end{figure}

Figure~\ref{fig:workflow} illustrates the overall workflow of the proposed GBNL framework for the prediction of protein-DNA/RNA binding. This approach couples algebraic invariants extracted from nucleic acid sequences with transformer-based representations of protein sequences, then learns a mapping from this combined representation to experimental binding affinities.

For each nucleotide, we construct position-based point clouds over a filtration and build Vietoris-Rips (VR) complexes with maximum simplex dimension 1. At each level, the complex is mapped to its Stanley-Reisner ideal and a minimal graded free resolution is computed to read off the graded Betti numbers \(\beta_{i,j}\). These \(\beta_{i,j}\) keys act as persistent algebraic signatures of local sequence organization and provide a compact, interpretable summary of multiscale combinatorial structure at the nucleotide level. In parallel, protein sequences are embedded with a modern transformer to capture global sequence regularities that complement the nucleotide specific algebraic features. The two feature streams are concatenated and used by a gradient boosted decision tree learning model to produce affinity predictions.

This design provides two benefits. First, the multiscale commutative algebra graded Betti numbers localize signal to specific nucleotides over filtration scales, which supports interpretation of model outputs. Second, the joint representation improves predictive accuracy by utilizing the rich transformer-based protein embeddings. Further details of protein sequence embedding, nucleic acid sequence processing, the datasets used and model parametrization are provided in the Methods section. 

\subsection{Protein-DNA/RNA Acid Binding Affinity Prediction}

\begin{figure}[t!]
	\centering
	\begin{subfigure}[b]{1\textwidth}
		\includegraphics[width=\textwidth]{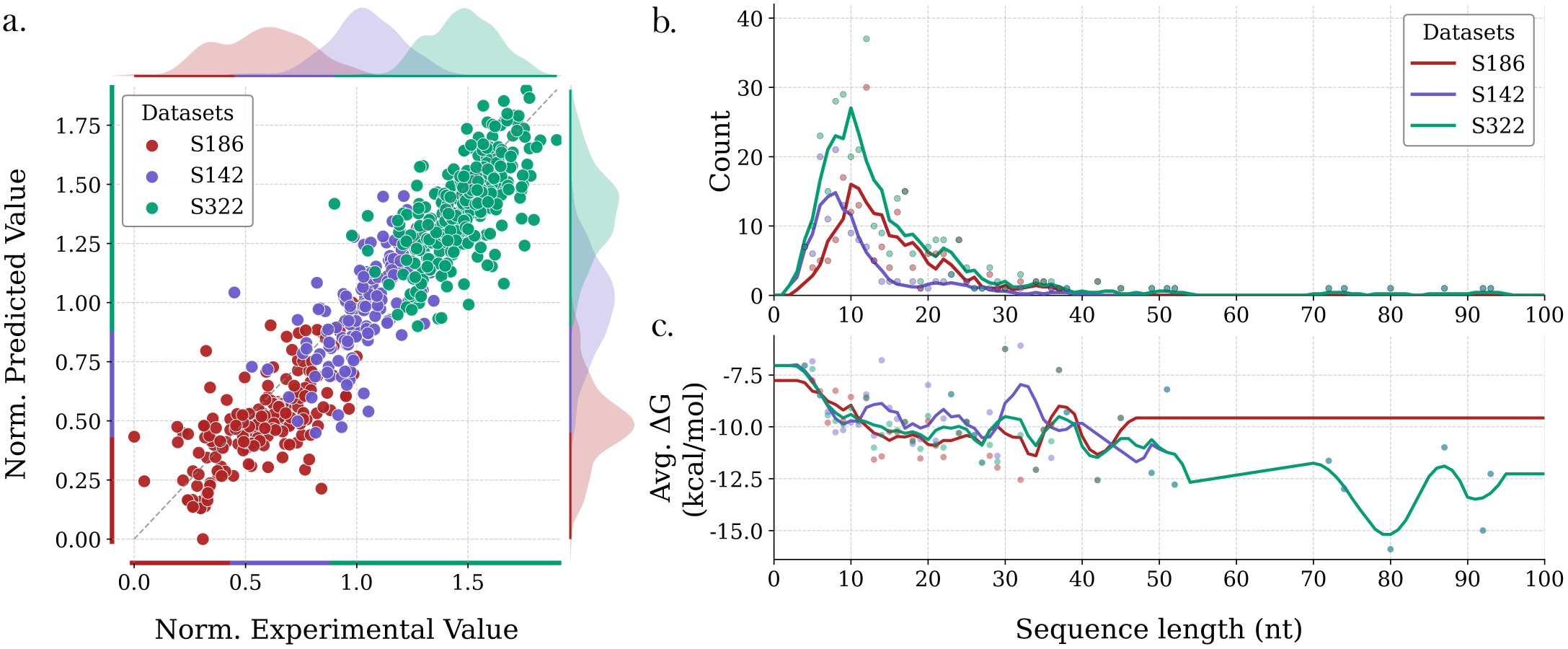}
	\end{subfigure}
	\caption{\textbf{Prediction summary and dataset statistics for S186, S142, and S322.}
		(a). Normalized experimental vs.\ predicted binding affinities; each dataset is min-max scaled to its own unit range and offset for visual separation, with marginal density plots shown along the top and right.
		(b). Sequence-length distribution: dots indicate the observed count at each length, and the curves show a 5-nt centered rolling average.
		(c). Experimental Gibbs free energy ($\Delta G$) vs.\ sequence length: dots mark the per-length mean $\Delta G$ values, and the curves show the corresponding 5-nt rolling average.}
	\label{fig:datasets}
\end{figure}

Protein-nucleic acid recognition drives catalysis, transport, signaling, transcription, and translation, while preserving chromosome stability and regulating gene expression. Disruption of these contacts is linked to autoimmune disease, inherited disorders, and cancer. Binding affinity reflects the combined effects of hydrogen bonding, van der Waals contacts, dipole interactions, electrostatics, and hydrophobic forces, and analyzing their joint impact enables rational, structure based therapeutic design. In this study, we use a graded Betti framework from commutative algebra to encode sequence derived determinants of these interactions and to predict protein-nucleic acid binding affinity with an interpretable, data driven model. Using this representation, we predict protein-DNA/RNA acid binding affinities with a model designed to balance accuracy with mechanistic insight.

We first assess predictive accuracy against the SVSBI reference~\cite{shen2023svsbi}. The SVSBI model integrates an ESM-based transformer for proteins with DNABERT for nucleic acids and reports an average Pearson of 0.669 and an RMSE of 1.98\,kcal/mol. On the other hand, our GBNL framework achieves \(R_p=0.691\) with RMSE \(=1.81\)\,kcal/mol on S186 using 3829 features, surpassing SVSBI on both correlation and error as shown in Table~\ref{table:Rp1}. Performance remains stable on S142 and S322, with 4287 and 4567 feature vectors, respectively, with \(R_p=0.66\), RMSE \(=2.17\)\,kcal/mol and \(R_p=0.657\), RMSE \(=2.03\)\,kcal/mol (Table~\ref{table:Rp2}). These results indicate consistent generalization across datasets of different sizes and compositions.

\autoref{fig:datasets}a compares experimental and predicted binding affinities after min-max normalization within each dataset, indicating stable calibration and no systematic bias between predictions and measurements. In \autoref{fig:datasets}b, sequence-length distributions are depicted with observed counts and smoothed 5-nt rolling averages. S186 exhibits a tight length profile centered near 12$\sim$nt with most sequences confined between 10-18$\sim$nt and dropping sharply beyond 20$\sim$nt. S142 is dominated by short fragments below 15$\sim$nt (median 9$\sim$nt) and showing fewer occurrences of longer strands. In contrast, S322 exhibits a broader mid-length representation, maintaining substantial counts across the 8-17$\sim$nt interval before gradually tapering at longer lengths, pointing to greater structural heterogeneity among complexes.

\begin{table}[htb!]
	\small
	\centering
	\caption{
		Comparison of prediction performance between existing SVSBI model and our GBNL framework on S186 for protein-nucleic acid binding affinity prediction. Reported metrics include Pearson correlation coefficient (\(R_p\)) and RMSE values in kcal/mol. All results are averaged over twenty independent runs with different random seeds, and the average metric values are reported.}
	\label{table:Rp1}
	\vspace{10pt}   
	\begin{tabular}{l | c | c}
		\hline
		\textbf{Model} & \textbf{\(R_p\)} & \textbf{RMSE (kcal/mol)} \\
		\hline
		\textbf{GBNL}   & \textbf{0.691} & \textbf{1.81} \\
		SVSBI \cite{shen2023svsbi}   & 0.669  &  1.98\footnotemark \\ 
		\hline
	\end{tabular}
\end{table}

\footnotetext{%
	The original RMSE of 1.45 reported by Shen et al. \cite{shen2023svsbi} was not converted into kcal/mol; here we apply the factor 1.3633 to obtain 1.98 kcal/mol.}

To illustrate how binding free energy (\(\Delta G\)) varies with sequence length across the datasets, \autoref{fig:datasets}c depicts the corresponding experimental and predicted trends. For S186, experimental energies span from \(-15.0\) to \(-5.8\)\,kcal/mol (mean \(-9.73 \pm 1.87\)\,kcal/mol), whereas predictions fall between \(-13.5\) and \(-5.8\)\,kcal/mol (mean \(-9.79 \pm 1.44\)\,kcal/mol). Both experimental and predicted profiles show a clear tendency toward more negative affinities as strand length increases, reaching minima near 42$\sim$nt around \(-12.6\)\,kcal/mol. A comparable pattern emerges in S142, where experimental \(\Delta G\) values range from \(-16.9\) to \(-4.3\)\,kcal/mol (mean \(-9.48 \pm 2.11\)\,kcal/mol), while predictions occupy a narrower interval of \(-13.3\) to \(-6.1\)\,kcal/mol (mean \(-9.48 \pm 1.46\)\,kcal/mol). Despite the strong short-length bias observed in its distribution, longer sequences exhibit more stable binding, with experimental and predicted minima near 80 and 93$\sim$nt at \(-15.9\) and \(-13.27\)\,kcal/mol, respectively. 

\begin{figure}[t!]
	\centering
	\begin{subfigure}[b]{1\textwidth}
		\includegraphics[width=\textwidth]{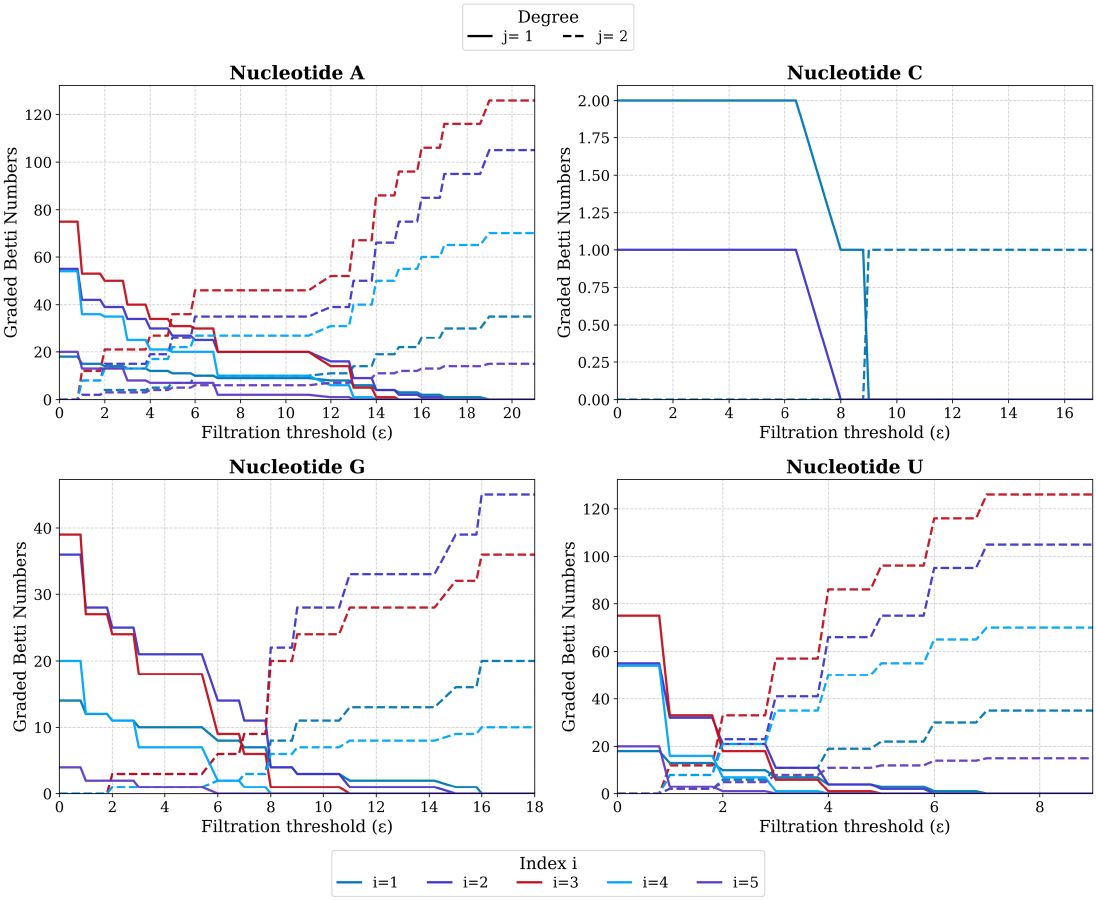}
	\end{subfigure}
	\caption{\textbf{Persistent graded Betti curves evolution across nucleotides for the RNA sequence \texttt{ACAUGUUUUCUGUGAAAACGGAG}.}
		Four panels correspond to nucleotides A, C, G, and U, showing how graded Betti counts vary with the filtration threshold~\(\varepsilon\). Curve groups are organized by homological index \(i\) and internal degree \(j\), as indicated in the legends.}
	\label{fig:ACGU}
\end{figure}

Furthermore, S322 follows the same overall trend, with experimental energies extending from \(-16.9\) to \(-4.3\)\,kcal/mol (mean \(-9.64 \pm 1.99\)\,kcal/mol) and predicted values spanning \(-13.4\) to \(-5.7\)\,kcal/mol (mean \(-9.62 \pm 1.46\)\,kcal/mol). Here too, both curves display increasingly negative \(\Delta G\) values with sequence length, culminating near 80–87$\sim$nt at roughly \(-13\) to \(-15.9\)\,kcal/mol. Collectively, these patterns indicate that stronger (more negative) experimental binding affinities are consistently associated with longer strands, even though the majority of sequences cluster at shorter lengths. Overall, the graded Betti framework aligns with the observed sequence-affinity relationships and maintains consistent calibration across datasets of different size and composition. On S186 in particular, our GBNL method delivers improved correlation and reduced error relative to the SVSBI baseline. A comprehensive comparison of experimental and predicted binding affinities for S186, S142, and S322 are presented in \autoref{tab:S186}, \autoref{tab:S142}, and \autoref{tab:S322}, respectively, along with the corresponding PDBIDs.

\begin{table}[htb!]
	\small
	\centering
	\caption{Prediction performance of GBNL on S142 and S322 datasets.}
	\label{table:Rp2}
	\vspace{8pt}
	\begin{tabular}{l | c | c}
		\hline
		\textbf{Dataset} & \textbf{\(R_p\)}  & \textbf{RMSE (kcal/mol)} \\
		\hline
		S142  & 0.66 & 2.17 \\
		S322 & 0.657 & 2.03 \\
		\hline
	\end{tabular}
\end{table}

To complement the quantitative evaluation and illustrate how graded features evolve across nucleotides, we visualize in \autoref{fig:ACGU} the persistent graded Betti curves in the GBNL framework for the RNA sequence \ttseq{ACAUGUUUUCUGUGAAAACGGAG} from S142. Each panel corresponds to a nucleotide (A, C, G, U) and shows how graded Betti counts vary with the filtration parameter~\(\varepsilon\), grouped by homological degree \(i\) and internal degree \(j\). Early curve separations and plateau differences reveal distinct algebraic organization among nucleotides, particularly in higher \((j-i)\) gaps, where broader dynamic ranges correspond to locally denser configurations. This visualization highlights the interpretive granularity of the graded Betti representation beyond aggregate performance metrics.

	\section{Discussion}\label{sec:Discussion}
	
	\subsection{Performance Analysis}
Building on these results, we find that Macaulay2 (M2) is generally efficient for graded Betti computations and scales well on higher performance cloud computing (HPCC) systems. However, performance degrades on longer or denser sequences because the associated Stanley-Reisner ideals develop very large generating sets and high syzygy complexity. An illustrative case in S186 is the protein-DNA complex 4LJR, which contains a 35$\sim$nt mononucleotide run of T. Under our position-based embedding, we mark the positions of a chosen token (e.g., all T’s) as 1-based indices along the sequence, and for each \(\varepsilon\) draw an edge between two positions if their index distance \(\le \varepsilon\). At \(\varepsilon=34\) these T positions become fully connected, so the VR(1) 1-skeleton is the complete graph \(K_{35}\) (here \(K_n\) denotes the complete graph on \(n\) vertices). For a graph-only complex, every triangle of \(K_{35}\) is a minimal nonface, hence a cubic generator of the Stanley-Reisner ideal; there are \(\binom{35}{3}=6545\) such generators in a single state. Resolving \(k[\Delta^{\epsilon}]=S/I^{\varepsilon}\) with thousands of degree-3 generators across many plateau endpoints in \(\varepsilon\) consumes time and memory. Computing maximum simplex dimension 2 is even more demanding because 4 or more cliques that are not promoted to 3-simplices become degree-4 and higher generators, accelerating the growth in the resolution.

To handle such cases while preserving the graded Betti output for graphs, we replace the computer algebra step with an exact combinatorial procedure specialized to VR(1). For each \(\varepsilon\), we set \(K=\lfloor \varepsilon \rfloor\) and work with the threshold graph that connects two indices when their gap is at most \(K\). Filtration events occur only at integer gaps, so Betti values are constant on plateaus and we evaluate at their right endpoints. In a VR(dim\(=1\)) complex, Hochster’s formula implies that only two diagonals of the graded Betti table can be nonzero, i.e., the disconnection diagonal \(\beta_{i,i+1}(k[\Delta(\varepsilon)])\), which counts how many \(j\)-vertex choices split into more than one component, and the cycle diagonal \(\beta_{i,i+2}(k[\Delta(\varepsilon)])\), which counts loops across all \(j\)-vertex choices,  additionally \(\beta_{0,0}=1\). We compute these exactly by direct counting: a small dynamic program yields the induced-subgraph component distribution for \(\beta_{i,i+1}\), and a closed identity combines total edge counts at threshold \(K\) with that distribution to obtain \(\beta_{i,i+2}\). Thus, the combinatorial evaluation reproduces the graded Betti tables for VR(1) while running time in seconds. We adopt this route only when M2 becomes prohibitively slow or explodes; for all other sequences we use M2, which remains efficient in our pipeline and benefits from HPCC using multiple jobs and subjobs.

We now place our GBNL findings in the context of prior sequence based predictors of protein-nucleic acid binding affinity. Existing models for protein-DNA \cite{zhao2010structure,yang2020predba,harini2023pda,yang2023empdba} and protein-RNA \cite{deng2019predprba,hong2023updated,harini2024pred} prediction have expanded the landscape. However, many depend on extensive hand-crafted descriptors, task specific tuning, and varied datasets, prohibiting direct comparison. Recent protein-DNA predictors range from potentials to ensembles. For example, emPDBA \cite{yang2023empdba} reports \(R_p=0.66\) with \(\mathrm{MAE}=1.24\,\mathrm{kcal/mol}\) after DNA subclassing (vs. \(R_p=0.12\), \(\mathrm{MAE}=1.64\,\mathrm{kcal/mol}\) unclassified). PDA-Pred \cite{harini2023pda} achieves \(R_p=0.78\) with categorization (vs. \(R_p=0.21\) without class categorization). On a 36-complex blind set, emPDBA attains \(R_p=0.53\) with \(\mathrm{MAE}=1.11\,\mathrm{kcal/mol}\); PreDBA \cite{yang2020predba} \(R_p=0.30\); DDNA3 \cite{zhao2010structure} \(R_p=0.09\) with \(\mathrm{MAE}=1.80\,\mathrm{kcal/mol}\). 

On the other hand, protein-RNA models are fewer. PRA-Pred \cite{harini2024pred} trains on 217 complexes clustered at 25\% identity, partitions RNA into five structural groups and proteins by function, and with jack-knife selection achieves \(R_p=0.77\) with \(\mathrm{MAE}=1.02\,\mathrm{kcal/mol}\); on the standard 44-complex blind set it reports \(R_p=0.60\) with \(\mathrm{MAE}=1.47\,\mathrm{kcal/mol}\). Similarly, PRdeltaGPred \cite{hong2023updated} reaches \(R_p=0.41\) with \(\mathrm{MAE}=1.83\,\mathrm{kcal/mol}\) while PredPRBA \cite{deng2019predprba} yields \(R_p=0.07\) with \(\mathrm{MAE}=2.07\,\mathrm{kcal/mol}\) (best subclass \(R_p=0.48\)). This demonstrates that comparison across literature is complicated by heterogeneous sequence identity thresholds and different subclass protocols utilized in each study. Moreover, for many methods, there are no reported codes, datasets, or protocol, making it difficult to reproduce the claimed results.  

In contrast, our GBNL framework uses curated, biochemically consistent datasets, including S186, S142 and a cleanly annotated extension, S322. It encodes paired sequences via algebraic grades and indices across multiscale filtrations and learns with a single gradient boosting regressor. Operating directly on primary sequences, the method requires no three-dimensional structures, subclass labels, or handcrafted interface features, yielding a single unified predictor for protein-DNA and protein-RNA binding affinity. Instead of enumerating hydrogen bonds, base-pair numbers, or stacking motifs, the GBNL framework summarizes interactions at the sequence level, avoiding reliance on subclass-specific models and enabling high-throughput analysis independent of curated annotations.

\begin{figure}[h!]
	\centering
	\begin{subfigure}[b]{1\textwidth}
		\includegraphics[width=\textwidth]{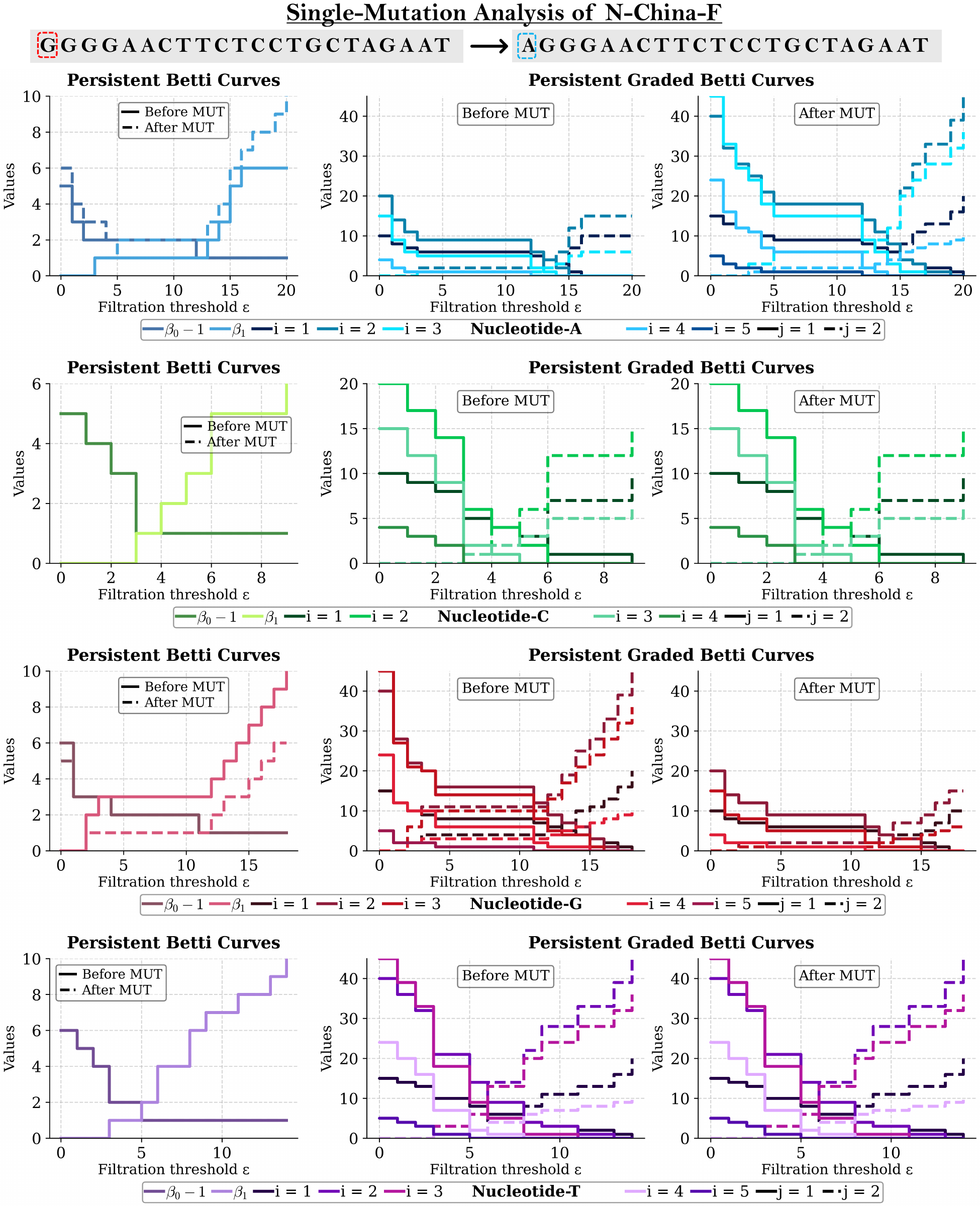}
	\end{subfigure}
	\caption{\textbf{Single mutation analysis by GBNL for N\text{-}China\text{-}F.}
		Rows: nucleotides. Column 1: PH Betti curves; columns 2 and 3: persistent graded Betti, all shown before and after mutation (MUT).}
	\label{fig:covid1}
\end{figure}

	\subsection{Model Interpretability}
	
Topological frameworks such as persistent homology (PH) is a well-established and powerful tool or multiscale quantification of biomolecular data. Yet, although PH captures global shifts in connectivity and cycles,it can obscure the underlying algebraic structures associated with specific mutations. Unlike PH summaries, GBNL separates homological content into graded components indexed by \(i\) and indices \(j\). Each \(\beta_{i,j}\) count generators (syzygies) at homological grade \(i\) and internal degree \(j\), yielding a precise account of how local perturbations shape the filtration landscape. Consequently, graded Betti curves uncovers intrinsic, nontrivial alterations that standard Betti summaries overlook.

We demonstrate in \autoref{fig:covid1} how GBNL sensitively captures a single-base mutation in the N-China-F primer of the SARS-CoV-2 N gene. A single nucleotide replacement in the reference sequence \ttseq{GGGGAACTTCTCCTGCTAGAAT} to \ttseq{AGGGAACTTCTCCTGCTAGAAT} at position 28881 (G\(\to\)A), one of the most frequent reported events (\(n=8672\)) \cite{wang2020mutations}, induces discernible shifts in topological and algebraic structure. At \(\varepsilon=0\), the A-specific component count increases from \(5\to6\) while loops remain zero. With larger \(\varepsilon\), PH shows modest rises in \(\beta_1\) (from \(2\to3\) at \(\varepsilon=13\) and \(6\to10\) at \(\varepsilon=20\)), indicating localized reorganization. Graded Betti profiles magnify these shifts: at the largest \(\varepsilon\), the gap \((j-i=2)\) for A increases \((1,3):10\to20\), \((2,4):15\to45\), \((3,5):6\to36\), while the corresponding G entries decrease \((1,3):20\to10\), \((2,4):45\to15\), \((3,5):36\to6\). Taken together, these opposing patterns reflect reallocation of generators and attenuation of higher-order syzygies near the mutated locus, aligning with a localized, single-site disturbance of the combinatorial neighborhood. PH and graded Betti curves for nucleotides C and T show no appreciable shift with mutation, confirming that all detected variation is attributable to the G\(\to\)A substitution. This analysis highlights the framework’s high sensitivity, showing that a single nucleotide alteration produces a measurable algebraic signature reflected in the graded Betti representation.

\begin{figure}[h!]
	\centering
	\begin{subfigure}[b]{1\textwidth}
		\includegraphics[width=\textwidth]{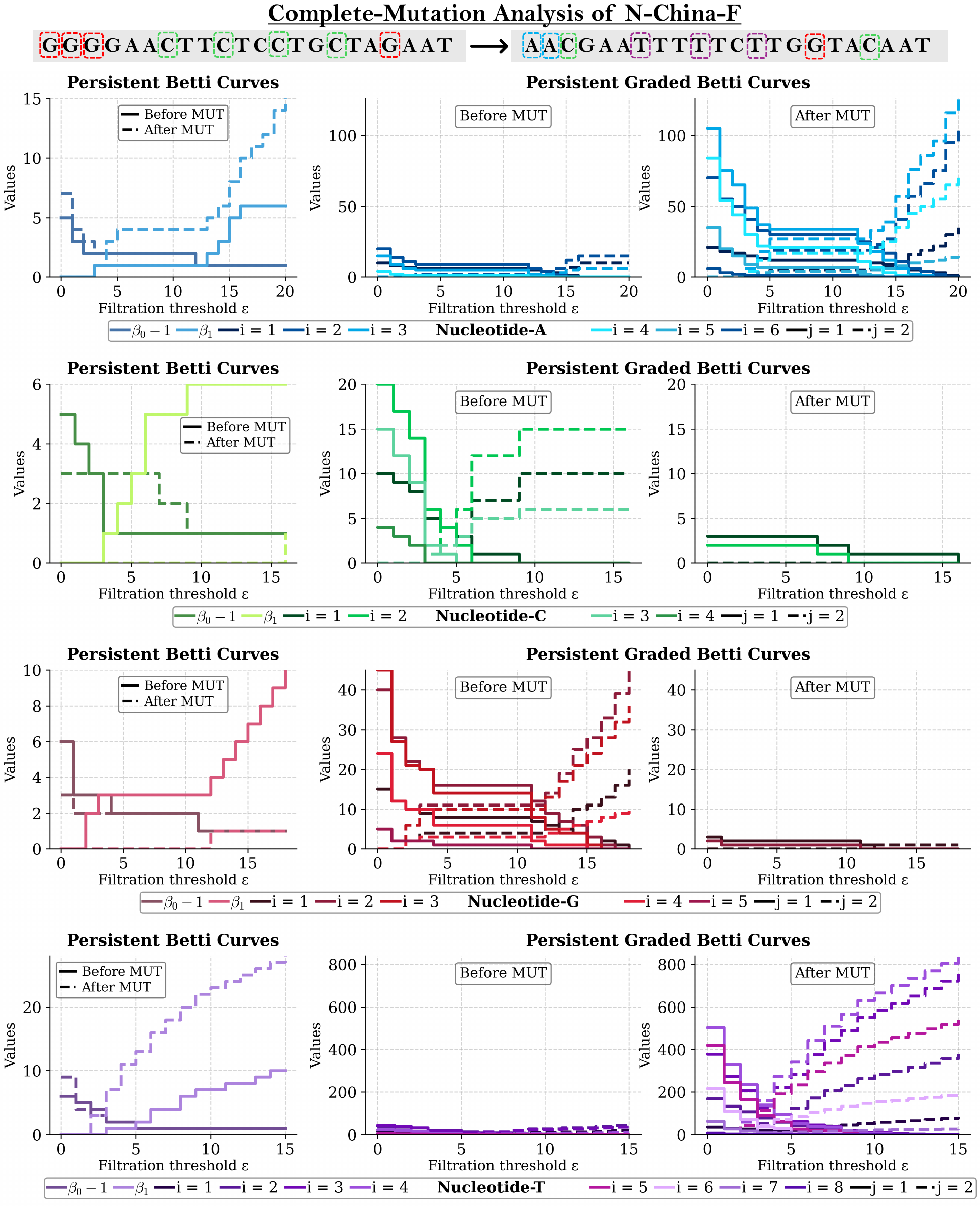}
	\end{subfigure}
	\caption{\textbf{Complete mutation analysis by GBNL for N\text{-}China\text{-}F.}
		Rows: nucleotides. Column 1: PH Betti curves; columns 2 and 3: persistent graded Betti, all shown before and after mutation (MUT).}
	\label{fig:covid_all}
\end{figure}

We next demonstrate how GBNL can efficiently capture the complete mutation in the same primer within the 22-base window. In \autoref{fig:covid_all}, the reference \ttseq{GGGGAACTTCTCCTGCTAGAAT} becomes \ttseq{AACGAATTTTTCTTGGTACAAT} \cite{wang2020mutations} after complete mutation. The triplet at 28881, 28882, and 28883 shows the highest counts (8672, 8660, 8659), and additional changes at 28887 (C\(\to\)T, 95), 28896 (C\(\to\)G, 52), 28899 (G\(\to\)C, 28), and 28900 (C\(\to\)T, 45) further perturb local connectivity. PH exhibits only limited variation in components and loops. Graded curves separate early in \(\varepsilon\), maintain broad plateaus, and span larger ranges at higher grades. For example, the component count shifts \(6\to9\) at \(\varepsilon=0\), and \(\beta_1\) increases from \(6\to20\) at \(\varepsilon=8\) and \(10\to28\) at \(\varepsilon=15\) for nucleotide T. At the largest \(\varepsilon\), the two-step gaps in the GBNL framework grow sharply: \((1,3):20\to84\), \((2,4):45\to378\), \((3,5):36\to756\), \((4,6):10\to840\). These values indicate a proliferation of algebraic generators across grades and scales. The graded Betti representation therefore reveals both precise single-site effects and the amplified structure induced by clustered substitutions. These cumulative changes illustrate how multiple concurrent mutations induce broad algebraic reorganization. The graded Betti representation exposing mutation-driven structural amplification far beyond what is visible through conventional persistent homology.

Together, the single- and complete-mutation analyses highlight the dual interpretive strengths of the GBNL framework: it resolves fine-grained single-site changes and it captures the algebraic amplification that arises from multiple, co-occurring substitutions within the same genomic window. This interpretability has direct modeling value. By encoding geometric scale through \(\varepsilon\) and algebraic organization through homological degree \(i\) and internal index \(j\), the graded Betti feature space links mutation-driven combinatorial rearrangements to functional outcomes. The clear, grade-resolved separations between reference and mutated trajectories in the graded columns of \autoref{fig:covid1} and \autoref{fig:covid_all} align with the observed gains in predictive accuracy when graded features are used instead of descriptors based only on persistent homology. Consequently, graded Betti numbers serve as mathematically rigorous commutative algebraic invariants and biologically interpretable indicators of mutation-driven reorganization in biomolecular sequences.

	\section{Methods}
	\label{sec:methods}
	
	In this section, we describe the datasets, the vectorization pipeline based on persistent commutative algebra, and transformer-based featurization via ESM2 protein embeddings. Next, we present the persistent Stanley-Reisner framework and show how it yields the $k$-mer algebraic representation of sequences. Finally, we summarize the machine learning models and evaluation metrics. 
	
	\subsection{Datasets} 
	
	In this study, we use three datasets. Dataset~S186, introduced by Shen et al.~\cite{shen2023svsbi}, contains 186 protein-nucleic acid complexes with strand lengths from 5 to 45$\sim$nt. It was assembled from PDBbind-v2020 under a strict curation protocol: a complex was retained only if it had one unique protein sequence and one unique nucleic acid sequence; multiple chains sharing the same sequence were allowed, but any chain with mixed or hybrid bases (for example, both T and U) was discarded. We kept only entries with experimentally measured binding affinities at 298~K, removed records with ambiguous labels (\textasciitilde, $<$, or $>$), and excluded complexes with strands shorter than five nucleotides. The second dataset, S142, is a protein-RNA collection created by merging our curated PDBbind-v2020 entries with selected PRA-Pred~\cite{harini2024pred} complexes and keeping only RNA-binding proteins; RNA strand lengths range from 4 to 93$\sim$nt. Finally, S322 is a combined protein-nucleic acid set constructed in the same manner, with strand lengths from 5 to 93$\sim$nt. Together, these curated data collections span from short to long strands and support the analysis of sequence-dependent binding affinities.
	
	\subsection{Persistent Commutative Algebra Featurization for Nucleic Acids}\label{sec:Vectorization}
	
	For each nucleic acid sequence, we convert per-nucleotide information into a fixed-length vector of graded Betti features using a uniform, discrete filtration. The occurrence positions of the four mononucleotides (A, C, G, T/U) are treated as one-dimensional point clouds and evaluated on an integer \(\varepsilon\)-grid. At each \(\varepsilon\), we build a Vietoris-Rips (VR) complex on the corresponding positions for a maximum simplicial dimension 1, form the associated Stanley-Reisner ideal, and compute the minimal free resolution to obtain graded Betti numbers \(\beta_{i,j}\). This yields, for every nucleotide, a sequence of \(\beta_{i,j}\) values aligned across identical filtration scales, providing uniformly sampled algebraic summaries of sequence organization.
	
	Each graded Betti index \((i,j)\) corresponds to a distinct position in the Betti table, characterized by its homological degree \(i\) and internal degree \(j\). Consistent feature alignment is implemented by constructing the feature vector from the global union of all nonzero \((i,j)\) indices observed over the \(\varepsilon\)-grid. By default, we exclude \(\beta_{0,0}\). We populate features using a left-endpoint carry-forward scheme: at each \(\varepsilon\), the entry records the nearest preceding (or equal) nonzero \(\beta_{i,j}\) to preserve consistent feature evolution across the \(\varepsilon\)-grid. For example, if for some nucleotide~C, we observe \(\beta_{1,2}=3\) and \(\beta_{2,3}=2\) at \(\varepsilon=0\), then \(\beta_{1,2}=2\) and \(\beta_{2,3}=1\) at \(\varepsilon=7\), and finally \(\beta_{1,2}=1\) at \(\varepsilon=9\), the entries corresponding to \(\beta_{1,2}\) across \(\varepsilon=0,\dots,9\) become \([3,3,3,3,3,3,3,2,2,1]\), while those for \(\beta_{2,3}\) become \([2,2,2,2,2,2,2,1,1,0]\). Stacking the nucleotide-specific blocks produces a unified representation of graded Betti changes across the filtration scales with an identical schema for every sequence.
	
	For S186, we observe that the global unions comprise 25, 29, 27, and 67 distinct \((i,j)\) keys for A, C, G, and T/U, summing to 148. Given ten filtration steps (\(\varepsilon = 0\ldots9\)), each sequence obtains a nucleic acid feature vector of length 148 × 10 = 1480. In S142, the per-nucleotide key counts are \(\{A:49,\,C:51,\,G:61,\,U:31\}\), yielding a 1920-dimensional RNA feature vector. In S322, the counts are \(\{A:49,\,C:51,\,G:61,\,T/U:67\}\), yielding a 2280-dimensional nucleic acid feature vector. These representations provide concise algebraic summaries of sequence structure while preserving key organizational features across datasets.
	
	\subsection{Transformer-Based Protein Language Model}\label{sec:NLP}
	
	Recent advances in natural language processing have introduced transformer architectures capable of capturing contextual dependencies directly from primary protein sequences. We employ the Evolutionary Scale Modeling (ESM2) framework~\cite{lin2023evolutionary}, a state-of-the-art protein language model trained on large-scale sequence corpora to learn residue-level contextual embeddings without structural supervision. ESM2 encodes each amino acid sequence through 36 transformer layers and produces 2560-dimensional embeddings that summarize long-range biochemical and evolutionary interactions.
	
	These protein embeddings complement the persistent commutative algebra descriptors of nucleic acids by providing a parallel, sequence-driven characterization of the interacting protein partner. For each protein-nucleic acid complex, the ESM2 protein embedding is concatenated with the corresponding nucleic-acid feature vector to form the final model input. The resulting combined feature dimensions are 3829 for S186, 4287 for S142, and 4567 for S322 after removing features that are identically zero across all samples. These final design matrices are employed as the inputs for model training and evaluation. \autoref{fig:workflow} provides a schematic overview of the feature-generation pipeline.
	
	\subsection{Persistent Stanley-Reisner Theory}
	\label{subsec:PSRT}
	
	Our framework leverages persistent Stanley-Reisner theory (PSRT) to ground data analysis in commutative algebra. Whereas traditional TDA uses persistent homology to reveal geometric and topological patterns like loops and voids~\cite{su2025topological}, PSRT instead examines the underlying algebraic and combinatorial essence of simplicial complexes. Input data are initially mapped onto a simplicial complex simplicial complex of vertices, edges, triangles, and higher-order simplices to retain their topological and combinatorial structure. A filtration is then applied to monitor the emergence and persistence of these characteristics across multiple spatial or geometric scales, producing commutative-algebraic invariants such as persistent $h$-vectors, $f$-vectors, graded Betti numbers, and facet ideals~\cite{suwayyid2025persistent}. This study concentrates solely on graded Betti numbers, situating the investigation within an algebra-driven framework for data analysis.
	
	\subsubsection{Persistent Stanley-Reisner Structures Over a Filtration}
	
	Let \(k\) be a field, and let \(\Delta\) be a simplicial complex on the finite vertex set \(V = \{x_1, \dots, x_n\}\). Suppose \(f: \Delta \to \mathbb{R}\) is a monotone function, i.e., \(f(\tau) \le f(\sigma)\) whenever \(\tau \subseteq \sigma\), which induces an increasing filtration \(\Delta^{s} \subseteq \Delta^{\epsilon}\) for \(s \le \varepsilon\). Let \(S = k[x_1, \dots, x_n]\) be the standard graded polynomial ring over \(k\), and for each \(\varepsilon \in \mathbb{R}\), define the {Stanley-Reisner ideal} of \(\Delta^{\epsilon}\) as
	\[
	I^{\varepsilon} := \left\langle x_{i_1} \cdots x_{i_r} \,\middle|\, \{x_{i_1}, \dots, x_{i_r}\} \notin \Delta^{\epsilon} \right\rangle \subseteq S,
	\]
	with corresponding {Stanley-Reisner ring}
	\[
	k[\Delta^{\epsilon}] := S / I^{\varepsilon}.
	\]
	
	Since the filtration is increasing, the subcomplexes satisfy \(\Delta^{s} \subseteq \Delta^{\epsilon}\) for \(s \le \varepsilon\), which implies a descending chain of monomial ideals:
	\[
	I^s \supseteq I^{\varepsilon} \quad \text{for all } s \le \varepsilon.
	\]

	\subsubsection{Persistent Graded Betti Numbers of Stanley-Reisner Rings}
	
	Let \(k\) be a field and \(S = k[x_1, \dots, x_n]\) the standard graded polynomial ring. For each filtration level \(\varepsilon \in \mathbb{R}\), the Stanley--Reisner ring \(k[\Delta^{\epsilon}] := S / I^{\varepsilon}\) inherits a natural \(\mathbb{Z}\)-graded \(S\)-module structure and admits a minimal graded free resolution:
	\begin{equation}\label{eq:min-free-res}
		\cdots \longrightarrow
		\bigoplus_{j} S(-j)^{\beta_{i,j}(k[\Delta^{\epsilon}])}
		\longrightarrow \cdots \longrightarrow
		k[\Delta^{\epsilon}] \longrightarrow 0,
	\end{equation}
	where \(\beta_{i,j}(k[\Delta^{\epsilon}]) := \dim_k \operatorname{Tor}^S_i(k[\Delta^{\epsilon}], k)_j\) are the {graded Betti numbers}.
	
	Hochster’s formula relates these graded Betti numbers to the topological Betti numbers of the induced subcomplexes:
	\begin{equation}\label{eq:Hochster-general}
		\beta_{i,j+i}(k[\Delta^{\epsilon}])
		=
		\sum_{\substack{W \subseteq V \\ |W| = j+i}} 
		\dim_k \widetilde{H}_{j-1}(\Delta_W^{\varepsilon}; k),
	\end{equation}
	where \(\widetilde{H}_{j-1}(\Delta_W^{\varepsilon}; k)\) denotes the \((j-1)\)-st reduced simplicial homology group over \(k\), and 
	\(\Delta_W^{\varepsilon} := \{ \sigma \in \Delta^{\varepsilon} \mid \sigma \subseteq W \}\) is the subcomplex induced on the vertex set \(W \subseteq V\).
	
	In particular, Hochster’s formula can be reformulated in terms of the (non-reduced) Betti numbers of induced subcomplexes. For each integer \(i \ge 0\), the following identities hold:
	\begin{align}
		\beta_{i,i+1}(k[\Delta^{\varepsilon}]) &= \sum_{\substack{W \subseteq V \\ |W| = i+1}} \left( \beta_0(\Delta_W^{\varepsilon}) - 1 \right), \label{eq:Hochster-j1} \\
		\beta_{i,i+j}(k[\Delta^{\varepsilon}]) &= \sum_{\substack{W \subseteq V \\ |W| = i+j}} \beta_{j-1}(\Delta_W^{\varepsilon}), \quad \text{for all } j \ge 2, \label{eq:Hochster-j2}
	\end{align}
	where \(\Delta_W^{\varepsilon}\) denotes the subcomplex of \(\Delta^{\varepsilon}\) induced on the vertex subset \(W \subseteq V\), and \(\beta_r(\Delta_W^{\varepsilon})\) denotes the \(r\)-th Betti number of \(\Delta_W^{\varepsilon}\) with coefficients in \(k\). 
	
	To refine this in a persistent setting, for \(\varepsilon \le \varepsilon'\), we define the {persistent graded Betti number}
	\begin{equation}\label{eq:persistent-Betti-short}
		\beta_{i, i+j}^{\varepsilon,\varepsilon'}(k[\Delta]) 
		:= \sum_{\substack{W \subseteq V \\ |W| = i + j}}
		\dim_k \left( \iota_{j-1}^{\varepsilon,\varepsilon'} : \widetilde{H}_{j-1}(\Delta_W^\varepsilon) 
		\to 
		\widetilde{H}_{j-1}(\Delta_W^{\varepsilon'}) \right),
	\end{equation}
	where \(\iota_{j-1}^{\varepsilon,\varepsilon'}\) is the homomorphism on reduced homology induced by inclusion. This provides a multigraded algebraic refinement of classical persistent Betti numbers, encoding both topological persistence and the combinatorial properties of the evolving homology classes.
	
	In the special case where \(|W| = |V|\), the persistent graded Betti number reduces to
	\[
	\beta_{i, |V|}^{\varepsilon,\varepsilon'} = \beta_{|V| - i - 1}^{\varepsilon,\varepsilon'},
	\]
	recovering the classical persistent Betti number of homological degree \(|V| - i - 1\). More generally, the family \(\{\beta_{i, i+j}^{\varepsilon,\varepsilon'}\}_{i,j}\) encodes a richer multiscale invariant that interpolates between algebraic and topological persistence.

	Additional structural identities among the persistent Betti numbers further simplify this formula. In particular, it is known that
	\[
	\beta^{\varepsilon,\varepsilon'}_{0,0} = 1, \qquad
	\beta^{\varepsilon,\varepsilon'}_{i,i} = 0 \quad \text{for all } i \ge 1, \qquad
	\beta^{\varepsilon,\varepsilon'}_{0,j} = 0 \quad \text{for all } j \ge 1, \qquad
	\beta^{\varepsilon,\varepsilon'}_{i,j} = 0 \quad \text{for all } i > j.
	\]
	Consequently, one obtains
	\[
	B_0 = \beta^{\varepsilon,\varepsilon'}_{0,0} = 1,
	\]
	and for each \( j \ge 1 \), the alternating sum simplifies to
	\[
	B_j = \sum_{i=1}^{j-1} (-1)^i \beta^{\varepsilon,\varepsilon'}_{i,j}.
	\]

	\subsubsection{\(k\)-mer Algebraic Representations of Sequences}
	\label{subsec:kmer_representation}
	
	In this section, we specialize the \(k\)-mer algebra framework to the case \(k=1\). The framework, introduced by Hozumi et al.~\cite{hozumi2024revealing}, provides a principled approach for embedding sequences as collections of integer sequences in a geometric space. Through the Stanley-Reisner construction, this embedding induces an algebraic structure determined by the positional distribution of individual nucleotides.
	
	Let \(\mathcal{A}\) be a finite alphabet. A \(1\)-mer over \(\mathcal{A}\) is an element \(\boldsymbol{x}=x_1\in\mathcal{A}\).
	Given a fixed \(1\)-mer \(\boldsymbol{x}\), define the indicator function
	\[
	\delta_{\boldsymbol{x}}:\mathcal{A}\to\{0,1\},
	\qquad
	\delta_{\boldsymbol{x}}(\boldsymbol{y})=
	\begin{cases}
		1,& \boldsymbol{y}=\boldsymbol{x},\\
		0,& \boldsymbol{y}\neq\boldsymbol{x}.
	\end{cases}
	\]
	For a sequence \(S=s_1s_2\cdots s_N\in\mathcal{A}^N\), the set of starting positions at which \(\boldsymbol{x}\) occurs is
	\[
	S^{\boldsymbol{x}}
	=
	\Bigl\{\, i\in[1,N] \ \Bigm|\ 
	\delta_{\boldsymbol{x}}\bigl(s_i\bigr)=1 \Bigr\}.
	\]
	We regard \(S^{\boldsymbol{x}}\subset\mathbb{R}\) as a one–dimensional point cloud and define the pairwise distance matrix
	\[
	D^{\boldsymbol{x}}=\bigl(d^{\boldsymbol{x}}_{ij}\bigr)_{i,j\in S^{\boldsymbol{x}}},
	\qquad
	d^{\boldsymbol{x}}_{ij}=|i-j|.
	\]
	
	These distance matrices serve as inputs to compute algebraic features over a filtration interval \([r_0,r_1]\subset\mathbb{R}_{\ge 0}\).
	For \(r,r'\in[r_0,r_1]\) with \(r\le r'\), build the Vietoris--Rips complex on \(S^{\boldsymbol{x}}\) and denote by
	\[
	v^{r,r'}_{\boldsymbol{x}}=\bigl( v^{r,r'}_{i,j}(\boldsymbol{x}) \bigr)_{i<j\in\mathbb{Z}_{\geq 0}}
	\]
	the vector of algebraic invariants in each dimension \(i\) and internal degree \(j\).
	On the diagonal \(r=r'\), write
	\[
	v_{\boldsymbol{x}}=\bigl( v_{i,j}(\boldsymbol{x}) \bigr)_{i<j\in\mathbb{Z}_{\geq 0}}
	\]
	to represent the curves for $r\in [r_0,r_1]$.
	The full \(1\)-mer representation of \(S\) is obtained by concatenation over all \(\boldsymbol{x}\in\mathcal{A}\):
	\[
	\boldsymbol{v}^1_S
	:=
	\bigl( v_{\boldsymbol{x}} \ \bigm|\ \boldsymbol{x}\in\mathcal{A} \bigr)
	=
	\bigl( v_{i,j}(\boldsymbol{x}) \ \bigm|\ \boldsymbol{x}\in\mathcal{A} \bigr)_{i<j\in\mathbb{Z}_{\geq 0}},
	\]
	which we refer to as the \(1\)-mer algebraic representation of \(S\).
	
	In this study, we fix the nucleotide alphabet 
	\(\mathcal{A} := \{A,\, C,\, G,\, T/U\}\).
	Accordingly, the \(1\)-mer algebraic representation of a sequence \(S\) is given by
	\[
	\boldsymbol{v}^1_S
	=
	\bigl(
	v_{\boldsymbol{A}},\,
	v_{\boldsymbol{C}},\,
	v_{\boldsymbol{G}},\,
	v_{\boldsymbol{T/U}}
	\bigr),
	\]
	where each component \(v_{\boldsymbol{x}}\) encodes the algebraic invariants associated with the positional distribution of the nucleotide \(\boldsymbol{x}\in\mathcal{A}\).

	\subsection{Examples}
	
	\begin{figure}[t!]
		\centering
		\begin{subfigure}[b]{0.24\textwidth}
			\includegraphics[width=\textwidth]{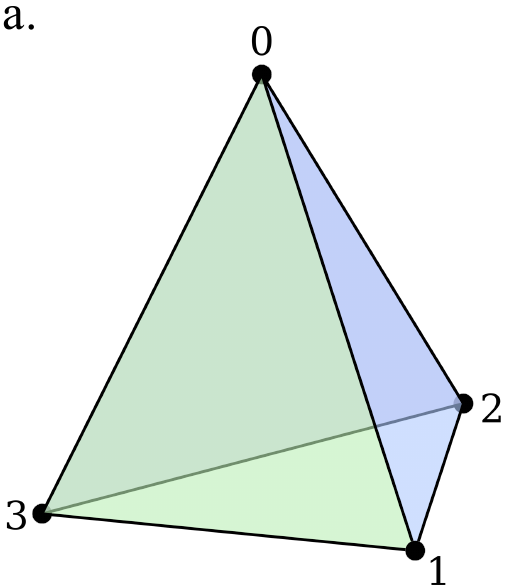}
		\end{subfigure}
		\hspace{6em} 
		\begin{subfigure}[b]{0.24\textwidth}
			\includegraphics[width=\textwidth]{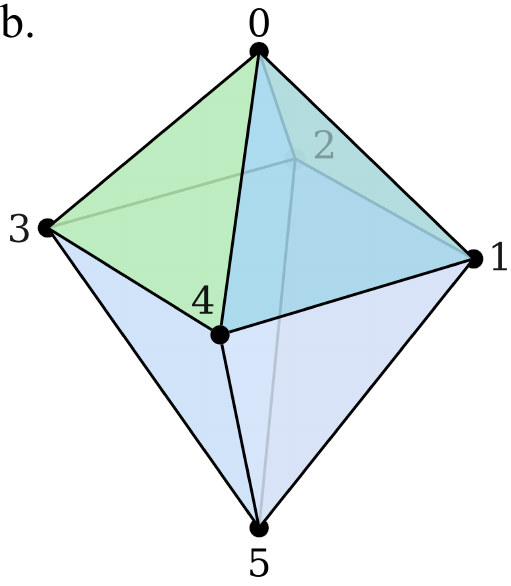}
		\end{subfigure}
		
		\caption{\textbf{Geometric realizations of the example complexes.}
			(a). Tetrahedron with one missing face. (b). Square bipyramid (octahedron boundary) with each base vertex connected to both apices.}
		
		\label{fig:examples}
	\end{figure}
	
	\begin{example}[Tetrahedron with one missing face]\label{example1}
		Let $\Delta$ be the simplicial complex on $V=\{0,1,2,3\}$ obtained from a tetrahedron by removing the base face $\{1,2,3\}$. 
		All six edges are present,
		\[
		\{0,1\},\ \{0,2\},\ \{0,3\},\ \{1,2\},\ \{1,3\},\ \{2,3\},
		\]
		and the $2$-faces (facets) are exactly
		\[
		\{0,1,2\},\ \{0,1,3\},\ \{0,2,3\}.
		\]

		Let $S=k[x_0,x_1,x_2,x_3]$ and write $I_\Delta\subset S$ for the Stanley--Reisner ideal.
		The only minimal nonface is $\{1,2,3\}$, so $I_\Delta=(x_1x_2x_3)$ is a principal (squarefree) cubic, and the minimal free resolution of $S/I_\Delta$ is
		\[
		0 \longrightarrow 
		\underbrace{S(-3)}_{i=1}
		\longrightarrow 
		\underbrace{S}_{i=0}
		\longrightarrow S/I_\Delta \longrightarrow 0.
		\]
		By the standard dictionary (each copy of $S(-j)$ in homological degree $i$ contributes one to $\beta_{i,j}$), we obtain
		\[
		\beta_{0,0}=1,\qquad \beta_{1,3}=1,
		\]
		and $\beta_{i,j}=0$ otherwise.
	\end{example}
	
	\begin{example}[Square bipyramid: topological vs.\ graded Betti numbers]\label{example2}
		Let $\Delta$ be the simplicial complex on the vertex set $V=\{0,1,2,3,4,5\}$ obtained as follows.
		Declare the \emph{square base} on $\{1,2,3,4\}$ with edges
		\[
		\{1,2\},\{2,3\},\{3,4\},\{4,1\},
		\]
		and connect each base vertex to both apices $0$ and $5$:
		\[
		\{0,i\},\{5,i\}\quad\text{for }i\in\{1,2,3,4\}.
		\]
		The $2$-faces (facets) are the eight triangles
		\[
		\{0,1,2\},\{0,2,3\},\{0,3,4\},\{0,4,1\},\qquad
		\{5,1,2\},\{5,2,3\},\{5,3,4\},\{5,4,1\}.
		\]
		Combinatorially, this is the boundary of the square bipyramid (the octahedron). 
		
		Hence over any field $k$,
		\[
		\tilde H_0(\Delta;k)=0,\quad \tilde H_1(\Delta;k)=0,\quad \tilde H_2(\Delta;k)\cong k,
		\]
		so the (unreduced) topological Betti numbers are
		\[
		\beta_0=1,\qquad \beta_1=0,\qquad \beta_2=1.
		\]
		Let $S=k[x_0,\dots,x_5]$ be standard graded and write $I_\Delta\subset S$ for the Stanley--Reisner ideal of $\Delta$.
		The minimal nonfaces of $\Delta$ are precisely the missing edges
		\[
		\{0,5\},\ \{1,3\},\ \{2,4\},
		\]
		whence
		\[
		I_\Delta=(x_0x_5,\ x_1x_3,\ x_2x_4).
		\]
		Because $I_\Delta$ is a complete intersection of three quadrics, the minimal free resolution of $S/I_\Delta$ over $S$ is the (pure) Koszul-type complex with shifts $2,2,2$:
		\[
		0 \longrightarrow 
		\underbrace{S(-6)}_{i=3}
		\longrightarrow 
		\underbrace{S(-4)^{3}}_{i=2}
		\longrightarrow 
		\underbrace{S(-2)^{3}}_{i=1}
		\longrightarrow 
		\underbrace{S}_{i=0}
		\longrightarrow S/I_\Delta \longrightarrow 0.
		\]
		Hence the graded Betti numbers are
		\[
		\beta_{0,0}=1,\quad \beta_{1,2}=3,\quad \beta_{2,4}=3,\quad \beta_{3,6}=1,
		\]
		and $\beta_{i,j}=0$ otherwise. Equivalently, the total Betti numbers are $(1,3,3,1)$.
	\end{example}

\subsection{Machine Learning Modeling}\label{sec:ML}

We train our regression models with 10-fold cross validation (CV) combined with Gradient Boosting Decision Trees (GBDT), implemented in Python via the \texttt{scikit-learn} library\,(v1.3.2). GBDT is valued for its resistance to overfitting, limited sensitivity to hyperparameter tuning, and straightforward deployment. It assembles many shallow decision trees generated from bootstrap samples of the training data and aggregates their outputs, so the ensemble corrects errors that any single tree might make. We supply the algorithm separately with PSRT-derived molecular descriptors and transformer-based descriptors. The hyperparameter settings are listed in \autoref{table:GBDT-parameters}.

\begin{table}[htb!]
	\small
	\centering
	\caption{Hyperparameters used in the gradient boosting regression model within a scikit-learn pipeline. StandardScaler is applied prior to model training.}
	\begin{tabular}{c c c c}
		\hline
		\textbf{No. of estimators} & \textbf{Max depth} & \textbf{Min. samples split} & \textbf{Learning rate} \\
		10,000 & 7 & 3 & 0.01 \\
		\hline
		\textbf{Max features} & \textbf{Subsample size} & \textbf{Random state} & \textbf{Standardization} \\
		Square root & 0.7 & Fixed (seeded) & Yes (StandardScaler) \\
		\hline
	\end{tabular}
	\label{table:GBDT-parameters}
\end{table}

\subsection{Evaluation Metrics}

We employ the Pearson correlation coefficient (\(R_p\)) for quantitative evaluation of the performance of our binding affinity prediction models, defined as: 
\begin{align*}
	\text{\(R_p\)}(\mathbf{x}, \mathbf{y}) = \frac{\sum_{m=1}^{M} (y_m^e - \bar{y}^e)(y_m^p - \bar{y}^p)}{\sqrt{\sum_{m=1}^{M} (y_m^e - \bar{y}^e)^2 \sum_{m=1}^{M} (y_m^p - \bar{y}^p)^2}},
\end{align*}
where \( y_m^e \) and \( y_m^p \) denote the experimental and predicted binding affinity values for the \( m \)-th sample, respectively, and \( \bar{y}^e \) and \( \bar{y}^p \) are their corresponding mean values.

We also report the root mean squared error (RMSE), which is computed as:
\begin{align*}
	\mathrm{RMSE} = \sqrt{\frac{1}{n} \sum_{m=1}^{M} (y_m^e - y_m^p)^2},
\end{align*}
where \( y_m^e \) and \( y_m^p \) represent the experimental and predicted binding affinity values for the \( m \)-th sample, respectively.

The above two metrics are employed to assess the performance of our machine learning models on all datasets. The original labels for these datasets are given as \( pK_d \) values, which are converted to binding free energies (in kcal/mol) by multiplying by a constant factor of 1.3633. Our models achieve reasonable RMSE values across all three datasets.

	\section{Conclusion}
	
	Commutative algebra provides a rigorous foundation for capturing combinatorial structure in biomolecular data. In this study, we introduce GBNL, a persistent commutative algebra framework for protein-DNA/RNA binding prediction.To our knowledge, this work presents the first graded-Betti application to large-data machine learning. DNA/RNA sequences are modeled as $k$-mer position sets, from which graded invariants of the associated Stanley-Reisner ideals are tracked over a positional filtration. In conjunction with transformer-based protein and nucleic acid embeddings, the algebraic features integrate nucleotide-level organization and protein-wide regularities into a concise, interpretable descriptor for predicting binding interactions.
	
	Beyond delivering precise affinity prediction, GBNL exhibits pronounced responsiveness to single-nucleotide mutations. For the N-China-F primer, a single nucleotide mutation yields a localized graded Betti shift that differentiates the mutated profile from the reference. A comprehensive mutation profiling of the same primer demonstrates that graded Betti signals capture subtle sequence perturbations beyond the resolution of standard Betti counts. Hence, GBNL upholds strong calibration and steady predictive behavior.. Furthermore, the GBNL achieves consistently strong performance by reaching a Pearson correlation of \(R_p = 0.69\) with RMSE \(= 1.81\)\,kcal/mol on S186, surpassing the SVSBI model~\cite{shen2023svsbi}, and maintaining comparable accuracy on S142 and S322. Collectively, these findings position GBNL as a broadly applicable and practical framework for biomolecular modeling.
	
	\section*{Data and Code availability}
	
	All codes needed to evaluate the conclusions in this study are available at \href{https://github.com/mzia-s/GBNL}{https://github.com/mzia-s/GBNL}. The three datasets can be found at \href{https://github.com/mzia-s/CAP}{https://github.com/mzia-s/CAP}. The SVSBI dataset can also be accessed at \href{https://github.com/WeilabMSU/SVS}{https://github.com/WeilabMSU/SVS}.
	
	\section*{Supporting Information}
	\href{./SI.pdf}{Supplementary Information} is available for supplementary tables. 
	
	\section*{Conflict of Interest}
	The authors declare no competing financial interests.

	\section*{Acknowledgments}
	This work was supported in part by NIH grant R35GM148196, NSF grant DMS-2052983, and MSU Research Foundation.    
	
	\clearpage

	\bibliographystyle{unsrt}
	\bibliography{refs}

\begin{thebibliography}{10}

\bibitem{bertoldo2023rna}
Jean~B Bertoldo, Simon M{\"u}ller, and Stefan H{\"u}ttelmaier.
\newblock Rna-binding proteins in cancer drug discovery.
\newblock {\em Drug discovery today}, 28(6):103580, 2023.

\bibitem{kisby2024introducing}
Glen~E Kisby, David~M Wilson~III, and Peter~S Spencer.
\newblock Introducing the role of genotoxicity in neurodegenerative diseases
  and neuropsychiatric disorders.
\newblock {\em International journal of molecular sciences}, 25(13):7221, 2024.

\bibitem{stockley2009filter}
Peter~G Stockley.
\newblock Filter-binding assays.
\newblock {\em DNA-Protein Interactions: Principles and Protocols, Third
  Edition}, pages 1--14, 2009.

\bibitem{zhang2005knowledge}
Chi Zhang, Song Liu, Qianqian Zhu, and Yaoqi Zhou.
\newblock A knowledge-based energy function for protein- ligand, protein-
  protein, and protein- dna complexes.
\newblock {\em Journal of medicinal chemistry}, 48(7):2325--2335, 2005.

\bibitem{nithin2019structure}
Chandran Nithin, Sunandan Mukherjee, and Ranjit~Prasad Bahadur.
\newblock A structure-based model for the prediction of protein--rna binding
  affinity.
\newblock {\em RNA}, 25(12):1628--1645, 2019.

\bibitem{bitencourt2019machine}
Gabriela Bitencourt-Ferreira and Walter~Filgueira de~Azevedo~Jr.
\newblock Machine learning to predict binding affinity.
\newblock In {\em Docking Screens for Drug Discovery}, pages 251--273.
  Springer, 2019.

\bibitem{hafner2010transcriptome}
Markus Hafner, Markus Landthaler, Lukas Burger, Mohsen Khorshid, Jean Hausser,
  Philipp Berninger, Andrea Rothballer, Manuel Ascano, Anna-Carina Jungkamp,
  Mathias Munschauer, et~al.
\newblock Transcriptome-wide identification of rna-binding protein and microrna
  target sites by par-clip.
\newblock {\em Cell}, 141(1):129--141, 2010.

\bibitem{zhao2011structure}
Huiying Zhao, Yuedong Yang, and Yaoqi Zhou.
\newblock Structure-based prediction of rna-binding domains and rna-binding
  sites and application to structural genomics targets.
\newblock {\em Nucleic acids research}, 39(8):3017--3025, 2011.

\bibitem{liu2016prediction}
Zhi-Ping Liu and Hongyu Miao.
\newblock Prediction of protein-rna interactions using sequence and structure
  descriptors.
\newblock {\em Neurocomputing}, 206:28--34, 2016.

\bibitem{shen2023svsbi}
Li~Shen, Hongsong Feng, Yuchi Qiu, and Guo-Wei Wei.
\newblock Svsbi: sequence-based virtual screening of biomolecular interactions.
\newblock {\em Communications biology}, 6(1):536, 2023.

\bibitem{tuszynska2011dars}
Irina Tuszynska and Janusz~M Bujnicki.
\newblock Dars-rnp and quasi-rnp: new statistical potentials for protein-rna
  docking.
\newblock {\em BMC bioinformatics}, 12:1--16, 2011.

\bibitem{setny2011coarse}
Piotr Setny and Martin Zacharias.
\newblock A coarse-grained force field for protein--rna docking.
\newblock {\em Nucleic acids research}, 39(21):9118--9129, 2011.

\bibitem{iwakiri2016improved}
Junichi Iwakiri, Michiaki Hamada, Kiyoshi Asai, and Tomoshi Kameda.
\newblock Improved accuracy in rna--protein rigid body docking by incorporating
  force field for molecular dynamics simulation into the scoring function.
\newblock {\em Journal of chemical theory and computation}, 12(9):4688--4697,
  2016.

\bibitem{lin2023evolutionary}
Zeming Lin, Halil Akin, Roshan Rao, Brian Hie, Zhongkai Zhu, Wenting Lu, Nikita
  Smetanin, Robert Verkuil, Ori Kabeli, Yaniv Shmueli, et~al.
\newblock Evolutionary-scale prediction of atomic-level protein structure with
  a language model.
\newblock {\em Science}, 379(6637):1123--1130, 2023.

\bibitem{song2024multiobjective}
Yao Song and Lu~Wang.
\newblock Multiobjective tree-based reinforcement learning for estimating
  tolerant dynamic treatment regimes.
\newblock {\em Biometrics}, 80(1):ujad017, 2024.

\bibitem{lo2018machine}
Yu-Chen Lo, Stefano~E Rensi, Wen Torng, and Russ~B Altman.
\newblock Machine learning in chemoinformatics and drug discovery.
\newblock {\em Drug discovery today}, 23(8):1538--1546, 2018.

\bibitem{cang2017topologynet}
Zixuan Cang and Guo-Wei Wei.
\newblock Topologynet: Topology based deep convolutional and multi-task neural
  networks for biomolecular property predictions.
\newblock {\em PLoS computational biology}, 13(7):e1005690, 2017.

\bibitem{papamarkou2024position}
Theodore Papamarkou, Tolga Birdal, Michael Bronstein, Gunnar Carlsson, Justin
  Curry, Yue Gao, Mustafa Hajij, Roland Kwitt, Pietro Lio, Paolo Di~Lorenzo,
  et~al.
\newblock Position: Topological deep learning is the new frontier for
  relational learning.
\newblock {\em Proceedings of machine learning research}, 235:39529, 2024.

\bibitem{nguyen2020review}
Duc~Duy Nguyen, Zixuan Cang, and Guo-Wei Wei.
\newblock A review of mathematical representations of biomolecular data.
\newblock {\em Physical Chemistry Chemical Physics}, 22(8):4343--4367, 2020.

\bibitem{nguyen2019mathematical}
Duc~Duy Nguyen, Zixuan Cang, Kedi Wu, Menglun Wang, Yin Cao, and Guo-Wei Wei.
\newblock Mathematical deep learning for pose and binding affinity prediction
  and ranking in d3r grand challenges.
\newblock {\em Journal of computer-aided molecular design}, 33:71--82, 2019.

\bibitem{nguyen2020mathdl}
Duc~Duy Nguyen, Kaifu Gao, Menglun Wang, and Guo-Wei Wei.
\newblock Mathdl: mathematical deep learning for d3r grand challenge 4.
\newblock {\em Journal of computer-aided molecular design}, 34(2):131--147,
  2020.

\bibitem{hozumi2024revealing}
Yuta Hozumi and Guo-Wei Wei.
\newblock Revealing the shape of genome space via k-mer topology.
\newblock {\em arXiv preprint arXiv:2412.20202}, 2024.

\bibitem{liu2025topological2}
Jian Liu, Li~Shen, Dong Chen, and Guo-Wei Wei.
\newblock Topological sequence analysis of genomes: Delta complex approaches.
\newblock {\em arXiv preprint arXiv:2507.05452}, 2025.

\bibitem{liu2025topological}
Jian Liu, Li~Shen, Mushal Zia, and Guo-Wei Wei.
\newblock Topological sequence analysis of genomes: Category approaches.
\newblock {\em arXiv preprint arXiv:2507.08043}, 2025.

\bibitem{zhao2010structure}
Huiying Zhao, Yuedong Yang, and Yaoqi Zhou.
\newblock Structure-based prediction of dna-binding proteins by structural
  alignment and a volume-fraction corrected dfire-based energy function.
\newblock {\em Bioinformatics}, 26(15):1857--1863, 2010.

\bibitem{rastogi2018accurate}
Chaitanya Rastogi, H~Tomas Rube, Judith~F Kribelbauer, Justin Crocker, Ryan~E
  Loker, Gabriella~D Martini, Oleg Laptenko, William~A Freed-Pastor, Carol
  Prives, David~L Stern, et~al.
\newblock Accurate and sensitive quantification of protein-dna binding
  affinity.
\newblock {\em Proceedings of the National Academy of Sciences},
  115(16):E3692--E3701, 2018.

\bibitem{barissi2022dnaffinity}
Sandro Barissi, Alba Sala, Mi{\l}osz Wiecz{\'o}r, Federica Battistini, and
  Modesto Orozco.
\newblock Dnaffinity: a machine-learning approach to predict dna binding
  affinities of transcription factors.
\newblock {\em Nucleic Acids Research}, 50(16):9105--9114, 2022.

\bibitem{yang2020predba}
Wenyi Yang and Lei Deng.
\newblock Predba: a heterogeneous ensemble approach for predicting protein-dna
  binding affinity.
\newblock {\em Scientific reports}, 10(1):1278, 2020.

\bibitem{harini2023pda}
K~Harini, Daisuke Kihara, and M~Michael Gromiha.
\newblock Pda-pred: Predicting the binding affinity of protein-dna complexes
  using machine learning techniques and structural features.
\newblock {\em Methods}, 213:10--17, 2023.

\bibitem{yang2023empdba}
Shuang Yang, Weikang Gong, Tong Zhou, Xiaohan Sun, Lei Chen, Wenxue Zhou, and
  Chunhua Li.
\newblock empdba: protein-dna binding affinity prediction by combining features
  from binding partners and interface learned with ensemble regression model.
\newblock {\em Briefings in Bioinformatics}, 24(4):bbad192, 2023.

\bibitem{yang2013dataset}
Xiufeng Yang, Haotian Li, Yangyu Huang, and Shiyong Liu.
\newblock The dataset for protein--rna binding affinity.
\newblock {\em Protein Science}, 22(12):1808--1811, 2013.

\bibitem{deng2019predprba}
Lei Deng, Wenyi Yang, and Hui Liu.
\newblock Predprba: prediction of protein-rna binding affinity using gradient
  boosted regression trees.
\newblock {\em Frontiers in genetics}, 10:637, 2019.

\bibitem{hong2023updated}
Xu~Hong, Xiaoxue Tong, Juan Xie, Pinyu Liu, Xudong Liu, Qi~Song, Sen Liu, and
  Shiyong Liu.
\newblock An updated dataset and a structure-based prediction model for
  protein--rna binding affinity.
\newblock {\em Proteins: Structure, Function, and Bioinformatics},
  91(9):1245--1253, 2023.

\bibitem{harini2024pred}
K~Harini, M~Sekijima, and M~Michael Gromiha.
\newblock Pra-pred: Structure-based prediction of protein-rna binding affinity.
\newblock {\em International Journal of Biological Macromolecules}, 259:129490,
  2024.

\bibitem{wee2025review}
JunJie Wee and Jian Jiang.
\newblock A review of topological data analysis and topological deep learning
  in molecular sciences.
\newblock {\em arXiv preprint arXiv:2509.16877}, 2025.

\bibitem{wang2020persistent}
Rui Wang, Duc~Duy Nguyen, and Guo-Wei Wei.
\newblock Persistent spectral graph.
\newblock {\em International journal for numerical methods in biomedical
  engineering}, 36(9):e3376, 2020.

\bibitem{chen2019evolutionary}
Jiahui Chen, Rundong Zhao, Yiying Tong, and Guo-Wei Wei.
\newblock Evolutionary de rham-hodge method.
\newblock {\em Discrete and continuous dynamical systems. Series B},
  26(7):3785, 2021.

\bibitem{memoli2022persistent}
Facundo M{\'e}moli, Zhengchao Wan, and Yusu Wang.
\newblock Persistent laplacians: Properties, algorithms and implications.
\newblock {\em SIAM Journal on Mathematics of Data Science}, 4(2):858--884,
  2022.

\bibitem{qiu2023persistent}
Yuchi Qiu and Guo-Wei Wei.
\newblock Persistent spectral theory-guided protein engineering.
\newblock {\em Nature computational science}, 3(2):149--163, 2023.

\bibitem{chen2022persistent}
Jiahui Chen, Yuchi Qiu, Rui Wang, and Guo-Wei Wei.
\newblock Persistent laplacian projected omicron ba. 4 and ba. 5 to become new
  dominating variants.
\newblock {\em Computers in Biology and Medicine}, 151:106262, 2022.

\bibitem{eisenbud2013commutative}
David Eisenbud.
\newblock {\em Commutative algebra: with a view toward algebraic geometry},
  volume 150.
\newblock Springer Science \& Business Media, 2013.

\bibitem{suwayyid2025persistent}
Faisal Suwayyid and Guo-Wei Wei.
\newblock Persistent stanley--reisner theory.
\newblock {\em arXiv preprint arXiv:2503.23482}, 2025.

\bibitem{feng2025caml}
Hongsong Feng, Faisal Suwayyid, Mushal Zia, JunJie Wee, Yuta Hozumi, Chun-Long
  Chen, and Guo-Wei Wei.
\newblock Caml: Commutative algebra machine learning---a case study on
  protein--ligand binding affinity prediction.
\newblock {\em Journal of Chemical Information and Modeling}, 2025.

\bibitem{wee2025commutative}
JunJie Wee, Faisal Suwayyid, Mushal Zia, Hongsong Feng, Yuta Hozumi, and
  Guo-Wei Wei.
\newblock Commutative algebra neural network reveals genetic origins of
  diseases.
\newblock {\em arXiv preprint arXiv:2509.26566}, 2025.

\bibitem{suwayyid2025cakl}
Faisal Suwayyid, Yuta Hozumi, Hongsong Feng, Mushal Zia, JunJie Wee, and
  Guo-Wei Wei.
\newblock Cakl: Commutative algebra k-mer learning of genomics.
\newblock {\em arXiv preprint arXiv:2508.09406}, 2025.

\bibitem{wang2020mutations}
Rui Wang, Yuta Hozumi, Changchuan Yin, and Guo-Wei Wei.
\newblock Mutations on covid-19 diagnostic targets.
\newblock {\em Genomics}, 112(6):5204--5213, 2020.

\bibitem{su2025topological}
Zhe Su, Xiang Liu, Layal~Bou Hamdan, Vasileios Maroulas, Jie Wu, Gunnar
  Carlsson, and Guo-Wei Wei.
\newblock Topological data analysis and topological deep learning beyond
  persistent homology--a review.
\newblock {\em arXiv preprint arXiv:2507.19504}, 2025.

\end{thebibliography}

\end{document}


\maketitle 
	\maketitle 
\newpage	
	{\setcounter{tocdepth}{4} \tableofcontents}
	\setcounter{page}{1}
	\newpage	
	
	\clearpage

	\section{Results: Experimental vs. Predicted Binding Affinities}
	\vspace{1em}
	
	\subsection{The S186 Dataset}

	\renewcommand*{\arraystretch}{1.2}
	\setlength{\tabcolsep}{6pt}
	
{\centering
\begin{longtable}{%
		>{\centering\arraybackslash}p{1.5cm}
		>{\centering\arraybackslash}p{1.5cm}
		>{\centering\arraybackslash}p{1.5cm}
		>{\centering\arraybackslash}p{1.5cm}
		>{\centering\arraybackslash}p{1.5cm}
		>{\centering\arraybackslash}p{1.5cm}
		>{\centering\arraybackslash}p{1.5cm}
		>{\centering\arraybackslash}p{1.5cm}
		>{\centering\arraybackslash}p{1.5cm}%
	}
	
\caption[Binding free energies for the PN dataset]%
{Experimental and predicted binding free energies ($\Delta G$) for the S186 dataset. \label{tab:S186}}\\

\toprule
\textbf{PDBID} & \textbf{Exp BA} & \textbf{Pred BA} &
\textbf{PDBID} & \textbf{Exp BA} & \textbf{Pred BA} &
\textbf{PDBID} & \textbf{Exp BA} & \textbf{Pred BA}\\
\midrule
\endfirsthead

\toprule
\textbf{PDBID} & \textbf{Exp BA} & \textbf{Pred BA} &
\textbf{PDBID} & \textbf{Exp BA} & \textbf{Pred BA} &
\textbf{PDBID} & \textbf{Exp BA} & \textbf{Pred BA}\\
\midrule
\endhead

1HVO & -6.760 & -7.225 & 1WET & -7.940 & -11.187 & 2PUF & -7.769 & -11.146 \\
1QPZ & -11.704 & -10.630 & 1QP0 & -10.165 & -10.900 & 1QP7 & -8.327 & -8.592 \\
1FJX & -12.070 & -12.670 & 1J75 & -10.165 & -9.320 & 1JFS & -11.803 & -10.718 \\
1J5K & -7.529 & -9.887 & 1PO6 & -9.276 & -9.719 & 1P51 & -11.528 & -11.275 \\
1QZG & -8.614 & -10.069 & 1OSB & -9.754 & -9.755 & 1HVN & -7.227 & -6.759 \\
1BDH & -8.242 & -8.067 & 1QP4 & -10.906 & -10.984 & 1QQB & -7.769 & -8.240 \\
1TW8 & -12.192 & -12.372 & 1QQA & -8.722 & -8.270 & 1DH3 & -12.070 & -10.208 \\
1JJ4 & -11.991 & -10.694 & 1FYM & -10.086 & -9.728 & 1JH9 & -11.859 & -10.335 \\
1JT0 & -9.965 & -10.364 & 1P78 & -11.528 & -11.543 & 1P71 & -11.528 & -11.450 \\
1OMH & -9.754 & -9.753 & 1S40 & -12.983 & -10.208 & 1U1L & -9.503 & -8.941 \\
1U1R & -9.053 & -9.092 & 1U1N & -8.862 & -9.708 & 1U1O & -9.855 & -9.680 \\
2BJC & -14.996 & -10.054 & 1ZZI & -7.713 & -9.822 & 2B0D & -8.242 & -9.112 \\
2GII & -12.783 & -12.783 & 2GIE & -12.592 & -12.394 & 2GIH & -12.252 & -12.256 \\
2ADY & -10.707 & -10.492 & 2AC0 & -10.364 & -10.676 & 2AYB & -11.859 & -11.752 \\
2ERE & -10.256 & -10.208 & 2GE5 & -8.751 & -8.998 & 2O6M & -10.751 & -9.906 \\
3D6Z & -8.204 & -6.799 & 2VYE & -9.991 & -10.179 & 3IGM & -8.590 & -9.099 \\
3HXQ & -12.572 & -12.576 & 3EQT & -9.514 & -9.685 & 3H15 & -7.455 & -9.078 \\
2KAE & -11.039 & -9.914 & 3MQ6 & -12.572 & -11.151 & 3JSP & -11.995 & -12.110 \\
3AAF & -8.994 & -9.450 & 3N1L & -11.407 & -11.236 & 3N1K & -11.236 & -11.236 \\
3R8F & -9.487 & -9.857 & 2RRA & -7.940 & -8.604 & 3QMI & -6.709 & -6.610 \\
2KXN & -7.700 & -8.934 & 1U1P & -9.107 & -9.095 & 1U1M & -8.901 & -9.128 \\
1U1K & -8.856 & -9.146 & 1U1Q & -9.706 & -9.321 & 2AOR & -8.242 & -8.332 \\
2AOQ & -7.438 & -8.160 & 2I9K & -12.983 & -11.835 & 2GIJ & -12.783 & -12.782 \\
2GIG & -12.252 & -12.253 & 2CCZ & -9.543 & -9.982 & 2AHI & -10.496 & -10.708 \\
2ERG & -10.364 & -10.108 & 2AYG & -11.859 & -11.862 & 2ES2 & -8.843 & -9.416 \\
2ATA & -9.646 & -10.484 & 2NP2 & -9.897 & -10.506 & 3D6Y & -6.799 & -8.203 \\
2VY1 & -9.573 & -10.163 & 3HXO & -12.572 & -12.581 & 3GIB & -12.070 & -9.539 \\
3D2W & -9.355 & -9.078 & 2KKF & -5.939 & -8.038 & 3M9E & -9.388 & -9.867 \\
3JSO & -12.402 & -11.580 & 3K3R & -11.966 & -12.335 & 3N1I & -11.554 & -11.385 \\
3N1J & -11.386 & -11.557 & 3KJP & -11.180 & -9.197 & 3Q0B & -8.134 & -9.051 \\
3PIH & -11.890 & -10.495 & 3QMH & -6.502 & -6.760 & 3QMB & -7.438 & -6.693 \\
3U7F & -8.482 & -9.038 & 3RN2 & -9.476 & -8.693 & 4ATK & -11.956 & -12.070 \\
2LTT & -10.628 & -10.133 & 4HIO & -10.467 & -10.133 & 4HIK & -10.388 & -10.100 \\
4A76 & -9.855 & -9.324 & 4HJ7 & -11.209 & -10.422 & 4HID & -8.273 & -10.433 \\
4LJ0 & -8.590 & -9.017 & 4HT8 & -9.471 & -11.510 & 4GCK & -9.256 & -9.693 \\
4F2J & -10.364 & -10.039 & 4GCT & -9.915 & -9.721 & 4NI7 & -13.223 & -10.045 \\
4CH1 & -9.003 & -9.251 & 4QJU & -8.768 & -10.865 & 4ZBN & -13.194 & -10.140 \\
4R55 & -7.876 & -9.646 & 4S0N & -10.751 & -9.369 & 4ZSF & -11.639 & -11.021 \\
4RKG & -6.258 & -8.954 & 5A72 & -9.916 & -10.088 & 5DWA & -12.332 & -10.690 \\
5K83 & -7.913 & -6.609 & 5ITH & -10.256 & -9.306 & 5W9S & -9.543 & -7.633 \\
5YI3 & -8.534 & -9.754 & 5XFP & -8.072 & -10.237 & 5MEY & -9.265 & -10.258 \\
5MEZ & -9.001 & -10.224 & 5W2M & -7.084 & -7.092 & 5K17 & -6.891 & -6.265 \\
6MG1 & -10.276 & -11.048 & 5VMV & -12.769 & -9.805 & 6CNQ & -8.242 & -8.518 \\
6FQP & -9.229 & -9.498 & 5WWF & -9.029 & -9.246 & 6FQQ & -8.739 & -9.538 \\
5ZVB & -5.990 & -5.807 & 5ZVA & -5.809 & -5.994 & 5MPF & -9.605 & -9.131 \\
6G1L & -12.332 & -11.597 & 6IIQ & -8.024 & -8.403 & 6IIR & -7.570 & -7.953 \\
3ON0 & -11.341 & -10.844 & 4HQU & -14.586 & -11.762 & 4HQX & -12.162 & -13.474 \\
4A75 & -9.232 & -9.782 & 3QSU & -10.511 & -10.320 & 4HJ8 & -10.440 & -8.443 \\
4HIM & -10.132 & -10.447 & 4HJ5 & -9.817 & -10.371 & 4HP1 & -7.028 & -8.847 \\
4J1J & -8.584 & -11.088 & 4NM6 & -8.140 & -9.217 & 4GCL & -9.278 & -9.781 \\
3ZPL & -11.751 & -10.058 & 3ZH2 & -10.057 & -10.890 & 4HT4 & -11.039 & -9.919 \\
4LNQ & -8.011 & -9.358 & 4LJR & -10.244 & -9.764 & 4R56 & -9.370 & -6.164 \\
4TMU & -12.030 & -7.977 & 4R22 & -10.819 & -10.276 & 4Z3C & -7.661 & -8.024 \\
3WPC & -10.496 & -11.423 & 3WPD & -11.619 & -10.488 & 2N8A & -9.573 & -9.975 \\
5HRT & -12.114 & -9.798 & 5T1J & -10.526 & -9.609 & 5VC9 & -9.265 & -7.914 \\
5HLG & -8.873 & -10.369 & 5W9Q & -8.482 & -7.551 & 6ASB & -7.940 & -9.285 \\
6ASD & -7.661 & -8.706 & 5K07 & -6.122 & -8.789 & 5YI2 & -9.738 & -8.680 \\
6MG3 & -11.449 & -10.443 & 6FWR & -8.873 & -10.264 & 6CNP & -8.391 & -8.704 \\
5ZD4 & -10.798 & -9.442 & 5ZMO & -9.163 & -9.928 & 6BWY & -6.251 & -9.687 \\
6CRM & -7.264 & -11.759 & 6CC8 & -7.181 & -9.422 & 6BUX & -5.807 & -7.535 \\
6A2I & -8.701 & -7.300 & 5ZKI & -8.015 & -9.841 & 6KBS & -7.713 & -9.830 \\
5ZKL & -10.479 & -10.538 & 5ZMD & -7.407 & -9.002 & 6ON0 & -9.605 & -9.870 \\
			
\bottomrule
\end{longtable}
}
	
	\newpage	
\subsection{The S142 Dataset}

\renewcommand*{\arraystretch}{1.2}
\setlength{\tabcolsep}{6pt}

{\centering
	\begin{longtable}{%
			>{\centering\arraybackslash}p{1.5cm}
			>{\centering\arraybackslash}p{1.5cm}
			>{\centering\arraybackslash}p{1.5cm}
			>{\centering\arraybackslash}p{1.5cm}
			>{\centering\arraybackslash}p{1.5cm}
			>{\centering\arraybackslash}p{1.5cm}
			>{\centering\arraybackslash}p{1.5cm}
			>{\centering\arraybackslash}p{1.5cm}
			>{\centering\arraybackslash}p{1.5cm}%
		}
	
\caption[Binding free energies for the PN dataset]%
{Experimental and predicted binding free energies ($\Delta G$) for the S142 dataset.\label{tab:S142}}\\
	
	\toprule
	\textbf{PDBID} & \textbf{Exp BA} & \textbf{Pred BA} &
	\textbf{PDBID} & \textbf{Exp BA} & \textbf{Pred BA} &
	\textbf{PDBID} & \textbf{Exp BA} & \textbf{Pred BA}\\
	\midrule
	\endfirsthead
	
	\toprule
	\textbf{PDBID} & \textbf{Exp BA} & \textbf{Pred BA} &
	\textbf{PDBID} & \textbf{Exp BA} & \textbf{Pred BA} &
	\textbf{PDBID} & \textbf{Exp BA} & \textbf{Pred BA}\\
	\midrule
	\endhead

2A9X & -9.303 & -10.473 & 3BX3 & -9.874 & -10.284 & 3IRW & -14.996 & -11.346 \\
3IWN & -12.270 & -13.266 & 3K0J & -10.992 & -12.732 & 3QGB & -11.528 & -11.272 \\
2XS2 & -10.113 & -7.344 & 3QGC & -10.982 & -11.091 & 4GHA & -8.590 & -9.518 \\
3V6Y & -11.464 & -11.558 & 2LUP & -6.090 & -9.969 & 4M59 & -9.531 & -12.078 \\
2MTV & -9.526 & -7.575 & 4RCM & -6.350 & -8.768 & 5EIM & -8.123 & -8.037 \\
5GXH & -7.661 & -9.352 & 5KLA & -10.334 & -11.677 & 5EN1 & -8.942 & -8.920 \\
5F5H & -9.668 & -9.548 & 5WZK & -9.265 & -10.278 & 5HO4 & -9.462 & -9.153 \\
5U9B & -9.388 & -8.399 & 5WZJ & -10.751 & -10.393 & 5TF6 & -11.639 & -10.197 \\
5SZE & -10.459 & -9.299 & 5WZG & -10.218 & -10.354 & 5M8I & -5.736 & -8.418 \\
5UDZ & -9.592 & -8.432 & 5YTV & -8.043 & -7.913 & 6FQ3 & -8.015 & -8.680 \\
5YTX & -8.007 & -8.025 & 5YKI & -12.030 & -9.498 & 6DCL & -10.647 & -9.177 \\
6GD3 & -5.191 & -8.436 & 6G2K & -6.427 & -8.578 & 6CMN & -11.727 & -9.590 \\
5YTS & -7.577 & -7.544 & 5WWG & -9.077 & -8.893 & 5WWF & -9.029 & -8.777 \\
5YTT & -7.394 & -7.730 & 6FQR & -6.073 & -7.331 & 6GX6 & -7.264 & -6.107 \\
5WWE & -8.550 & -9.052 & 6GC5 & -8.435 & -6.123 & 6RT6 & -6.755 & -7.439 \\
6NOF & -10.617 & -10.476 & 6R7B & -10.558 & -9.283 & 6NOC & -10.490 & -10.608 \\
6NOH & -10.421 & -10.598 & 6A6J & -9.068 & -8.604 & 6NOD & -10.280 & -11.411 \\
6G99 & -6.176 & -8.955 & 6RT7 & -7.405 & -7.235 & 6U9X & -10.086 & -9.517 \\
6NY5 & -9.952 & -10.058 & 6GBM & -5.556 & -9.276 & 1EC6 & -7.990 & -8.215 \\
1M8Y & -9.270 & -12.626 & 1RPU & -13.300 & -9.854 & 1UTD & -16.890 & -9.025 \\
2B6G & -10.600 & -9.883 & 2ERR & -12.200 & -8.730 & 2F8K & -10.300 & -9.521 \\
2G4B & -7.830 & -8.004 & 2KFY & -8.840 & -8.235 & 2KG0 & -7.400 & -8.317 \\
2KG1 & -7.450 & -8.066 & 2KX5 & -11.400 & -10.711 & 2KXN & -8.170 & -8.520 \\
2L41 & -4.250 & -7.391 & 2LA5 & -11.500 & -10.650 & 2LEB & -9.320 & -9.403 \\
2LEC & -9.440 & -9.346 & 2M8D & -8.730 & -8.557 & 2MJH & -9.670 & -8.932 \\
2MXY & -7.900 & -6.865 & 2MZ1 & -6.870 & -7.919 & 2N82 & -10.300 & -9.990 \\
2RRA & -7.940 & -8.340 & 2RU3 & -8.520 & -9.249 & 2XC7 & -7.040 & -8.415 \\
2XFM & -8.240 & -7.512 & 2XNR & -5.400 & -7.238 & 2ZKO & -8.050 & -9.846 \\
3BSB & -10.200 & -11.009 & 3BSX & -11.500 & -10.532 & 3BX2 & -9.970 & -9.886 \\
3EQT & -9.520 & -9.562 & 3GIB & -10.690 & -9.069 & 3K49 & -12.050 & -10.890 \\
3K4E & -10.800 & -11.583 & 3K5Q & -9.300 & -8.987 & 3K5Y & -10.540 & -9.138 \\
3K5Z & -9.380 & -10.044 & 3K61 & -8.890 & -9.736 & 3K62 & -8.740 & -9.416 \\
3K64 & -9.220 & -9.080 & 3L25 & -8.175 & -9.889 & 3LQX & -12.220 & -11.041 \\
3MDG & -7.850 & -8.541 & 3MOJ & -13.000 & -10.866 & 3NCU & -10.300 & -9.436 \\
3NMR & -8.440 & -7.327 & 3NNH & -6.190 & -8.412 & 3O3I & -6.520 & -7.071 \\
3O6E & -7.070 & -6.522 & 3Q0L & -12.300 & -11.563 & 3Q0M & -12.500 & -12.203 \\
3Q0N & -10.500 & -11.819 & 3Q0P & -12.800 & -10.550 & 3Q0Q & -13.300 & -12.919 \\
3Q0R & -13.800 & -12.187 & 3Q0S & -10.800 & -13.100 & 3QG9 & -10.600 & -11.410 \\
3U4M & -15.900 & -11.271 & 3V71 & -11.100 & -10.082 & 3V74 & -11.600 & -11.444 \\
3WBM & -9.980 & -10.359 & 4CIO & -9.760 & -8.356 & 4ED5 & -9.080 & -8.419 \\
4ERD & -9.540 & -9.431 & 4HT8 & -9.380 & -10.120 & 4JVH & -9.590 & -8.857 \\
4KJI & -8.940 & -10.143 & 4LG2 & -9.540 & -8.546 & 4NL3 & -7.380 & -10.311 \\
4O26 & -8.200 & -9.945 & 4OE1 & -12.110 & -9.528 & 4QI2 & -9.440 & -9.448 \\
4QVC & -11.560 & -10.141 & 4QVD & -10.120 & -11.080 & 4R3I & -7.770 & -7.439 \\
4RCJ & -8.180 & -8.683 & 4TUX & -10.400 & -9.879 & 4U8T & -9.130 & -7.510 \\
4Z31 & -8.790 & -9.128 & 5DNO & -7.830 & -8.487 & 5ELR & -6.250 & -9.016 \\
5V7C & -6.990 & -9.281 & 5W1I & -12.790 & -10.724 & 5WZH & -10.460 & -10.501 \\
6D12 & -9.300 & -10.253 & & &  & & &  \\
	
\bottomrule
\end{longtable}
}
	
\newpage	
\subsection{The S322 Dataset}

\renewcommand*{\arraystretch}{1.2}
\setlength{\tabcolsep}{6pt}

\centering
\begin{longtable}{%
		>{\centering\arraybackslash}p{1.5cm}
		>{\centering\arraybackslash}p{1.5cm}
		>{\centering\arraybackslash}p{1.5cm}
		>{\centering\arraybackslash}p{1.5cm}
		>{\centering\arraybackslash}p{1.5cm}
		>{\centering\arraybackslash}p{1.5cm}
		>{\centering\arraybackslash}p{1.5cm}
		>{\centering\arraybackslash}p{1.5cm}
		>{\centering\arraybackslash}p{1.5cm}}
	
\caption[Binding free energies for the PN dataset]%
{Experimental and predicted binding free energies ($\Delta G$) for the S322 dataset.\label{tab:S322}}\\
	
	\toprule
	\textbf{PDBID} & \textbf{Exp BA} & \textbf{Pred BA} &
	\textbf{PDBID} & \textbf{Exp BA} & \textbf{Pred BA} &
	\textbf{PDBID} & \textbf{Exp BA} & \textbf{Pred BA}\\
	\midrule
	\endfirsthead
	
	\toprule
	\textbf{PDBID} & \textbf{Exp BA} & \textbf{Pred BA} &
	\textbf{PDBID} & \textbf{Exp BA} & \textbf{Pred BA} &
	\textbf{PDBID} & \textbf{Exp BA} & \textbf{Pred BA}\\
	\midrule
	\endhead
	
1HVO & -6.760 & -7.226 & 1WET & -7.940 & -11.440 & 2PUF & -7.769 & -11.062 \\
1QPZ & -11.704 & -10.478 & 1QP0 & -10.165 & -11.059 & 1QP7 & -8.327 & -8.090 \\
1FJX & -12.070 & -12.625 & 1J75 & -10.165 & -8.913 & 1JFS & -11.803 & -10.968 \\
1J5K & -7.529 & -8.077 & 1PO6 & -9.276 & -9.632 & 1P51 & -11.528 & -11.447 \\
1QZG & -8.614 & -9.820 & 1OSB & -9.754 & -9.755 & 1HVN & -7.227 & -6.768 \\
1BDH & -8.242 & -8.273 & 1QP4 & -10.906 & -10.964 & 1QQB & -7.769 & -8.822 \\
1TW8 & -12.192 & -12.403 & 1QQA & -8.722 & -7.864 & 1DH3 & -12.070 & -10.220 \\
1JJ4 & -11.991 & -10.522 & 1FYM & -10.086 & -9.240 & 1JH9 & -11.859 & -8.205 \\
1JT0 & -9.965 & -10.360 & 1P78 & -11.528 & -11.575 & 1P71 & -11.528 & -11.430 \\
1OMH & -9.754 & -9.754 & 1S40 & -12.983 & -9.734 & 1U1L & -9.503 & -8.950 \\
1U1R & -9.053 & -9.083 & 1U1N & -8.862 & -9.706 & 1U1O & -9.855 & -9.703 \\
2BJC & -14.996 & -9.816 & 1ZZI & -7.713 & -7.881 & 2B0D & -8.242 & -9.204 \\
2GII & -12.783 & -12.769 & 2GIE & -12.592 & -12.424 & 2GIH & -12.252 & -12.279 \\
2ADY & -10.707 & -10.470 & 2AC0 & -10.364 & -10.593 & 2AYB & -11.859 & -11.559 \\
2ERE & -10.256 & -10.138 & 2GE5 & -8.751 & -9.085 & 2O6M & -10.751 & -9.741 \\
3D6Z & -8.204 & -6.800 & 2VYE & -9.991 & -9.808 & 3IGM & -8.590 & -8.929 \\
3HXQ & -12.572 & -12.333 & 3EQT & -9.514 & -9.508 & 3H15 & -7.455 & -8.929 \\
2KAE & -11.039 & -9.052 & 3MQ6 & -12.572 & -11.236 & 3JSP & -11.995 & -11.966 \\
3AAF & -8.994 & -9.124 & 3N1L & -11.407 & -11.237 & 3N1K & -11.236 & -11.407 \\
3R8F & -9.487 & -10.585 & 2RRA & -7.940 & -8.340 & 3QMI & -6.709 & -6.576 \\
2KXN & -7.700 & -8.578 & 1U1P & -9.107 & -9.067 & 1U1M & -8.901 & -8.950 \\
1U1K & -8.856 & -9.146 & 1U1Q & -9.706 & -9.325 & 2AOR & -8.242 & -7.950 \\
2AOQ & -7.438 & -8.185 & 2I9K & -12.983 & -11.701 & 2GIJ & -12.783 & -12.756 \\
2GIG & -12.252 & -12.261 & 2CCZ & -9.543 & -9.670 & 2AHI & -10.496 & -10.680 \\
2ERG & -10.364 & -10.131 & 2AYG & -11.859 & -11.871 & 2ES2 & -8.843 & -8.864 \\
2ATA & -9.646 & -10.486 & 2NP2 & -9.897 & -10.332 & 3D6Y & -6.799 & -8.203 \\
2VY1 & -9.573 & -9.619 & 3HXO & -12.572 & -10.434 & 3GIB & -12.070 & -9.425 \\
3D2W & -9.355 & -8.362 & 2KKF & -5.939 & -8.160 & 3M9E & -9.388 & -9.679 \\
3JSO & -12.402 & -11.844 & 3K3R & -11.966 & -12.410 & 3N1I & -11.554 & -11.391 \\
3N1J & -11.386 & -11.543 & 3KJP & -11.180 & -9.107 & 3Q0B & -8.134 & -9.133 \\
3PIH & -11.890 & -10.726 & 3QMH & -6.502 & -6.783 & 3QMB & -7.438 & -6.892 \\
3U7F & -8.482 & -8.521 & 3RN2 & -9.476 & -8.758 & 4ATK & -11.956 & -12.002 \\
2LTT & -10.628 & -10.323 & 4HIO & -10.467 & -10.180 & 4HIK & -10.388 & -10.364 \\
4A76 & -9.855 & -9.274 & 4HJ7 & -11.209 & -10.408 & 4HID & -8.273 & -10.425 \\
4LJ0 & -8.590 & -9.068 & 4HT8 & -9.471 & -10.284 & 4GCK & -9.256 & -9.689 \\
4F2J & -10.364 & -9.878 & 4GCT & -9.915 & -9.858 & 4NI7 & -13.223 & -10.345 \\
4CH1 & -9.003 & -8.548 & 4QJU & -8.768 & -11.131 & 4ZBN & -13.194 & -10.103 \\
4R55 & -7.876 & -9.325 & 4S0N & -10.751 & -9.301 & 4ZSF & -11.639 & -11.165 \\
4RKG & -6.258 & -8.980 & 5A72 & -9.916 & -10.278 & 5DWA & -12.332 & -10.896 \\
5K83 & -7.913 & -7.264 & 5ITH & -10.256 & -8.855 & 5W9S & -9.543 & -7.710 \\
5YI3 & -8.534 & -9.698 & 5XFP & -8.072 & -9.177 & 5MEY & -9.265 & -9.167 \\
5MEZ & -9.001 & -9.384 & 5W2M & -7.084 & -6.741 & 5K17 & -6.891 & -8.259 \\
6MG1 & -10.276 & -10.966 & 5VMV & -12.769 & -9.771 & 6CNQ & -8.242 & -8.108 \\
6FQP & -9.229 & -9.097 & 5WWF & -9.029 & -8.842 & 6FQQ & -8.739 & -9.350 \\
5ZVB & -5.990 & -5.745 & 5ZVA & -5.809 & -6.015 & 5MPF & -9.605 & -8.895 \\
6G1L & -12.332 & -11.340 & 6IIQ & -8.024 & -8.503 & 6IIR & -7.570 & -8.067 \\
3ON0 & -11.341 & -10.682 & 4HQU & -14.586 & -11.618 & 4HQX & -12.162 & -13.372 \\
4A75 & -9.232 & -9.764 & 3QSU & -10.511 & -9.271 & 4HJ8 & -10.440 & -8.476 \\
4HIM & -10.132 & -10.458 & 4HJ5 & -9.817 & -10.388 & 4HP1 & -7.028 & -8.723 \\
4J1J & -8.584 & -10.830 & 4NM6 & -8.140 & -8.601 & 4GCL & -9.278 & -9.690 \\
3ZPL & -11.751 & -10.363 & 3ZH2 & -10.057 & -10.904 & 4HT4 & -11.039 & -10.697 \\
4LNQ & -8.011 & -9.399 & 4LJR & -10.244 & -9.695 & 4R56 & -9.370 & -7.009 \\
4TMU & -12.030 & -8.235 & 4R22 & -10.819 & -10.117 & 4Z3C & -7.661 & -8.099 \\
3WPC & -10.496 & -11.349 & 3WPD & -11.619 & -10.537 & 2N8A & -9.573 & -10.712 \\
5HRT & -12.114 & -9.735 & 5T1J & -10.526 & -9.740 & 5VC9 & -9.265 & -7.629 \\
5HLG & -8.873 & -10.472 & 5W9Q & -8.482 & -7.690 & 6ASB & -7.940 & -8.815 \\
6ASD & -7.661 & -8.524 & 5K07 & -6.122 & -7.002 & 5YI2 & -9.738 & -8.693 \\
6MG3 & -11.449 & -10.123 & 6FWR & -8.873 & -9.624 & 6CNP & -8.391 & -8.247 \\
5ZD4 & -10.798 & -9.382 & 5ZMO & -9.163 & -9.646 & 6BWY & -6.251 & -9.902 \\
6CRM & -7.264 & -12.002 & 6CC8 & -7.181 & -8.992 & 6BUX & -5.807 & -7.637 \\
6A2I & -8.701 & -7.286 & 5ZKI & -8.015 & -9.564 & 6KBS & -7.713 & -9.610 \\
5ZKL & -10.479 & -10.343 & 5ZMD & -7.407 & -9.140 & 6ON0 & -9.605 & -10.027 \\
2A9X & -9.303 & -10.690 & 3BX3 & -9.874 & -10.327 & 3IRW & -14.996 & -10.938 \\
3IWN & -12.270 & -10.902 & 3K0J & -10.992 & -13.127 & 3QGB & -11.528 & -11.195 \\
2XS2 & -10.113 & -8.308 & 3QGC & -10.982 & -11.351 & 4GHA & -8.590 & -9.688 \\
3V6Y & -11.464 & -11.553 & 2LUP & -6.090 & -10.774 & 4M59 & -9.531 & -12.018 \\
2MTV & -9.526 & -7.708 & 4RCM & -6.350 & -9.235 & 5EIM & -8.123 & -8.468 \\
5GXH & -7.661 & -9.430 & 5KLA & -10.334 & -11.633 & 5EN1 & -8.942 & -9.055 \\
5F5H & -9.668 & -9.147 & 5WZK & -9.265 & -10.225 & 5HO4 & -9.462 & -9.166 \\
5U9B & -9.388 & -8.906 & 5WZJ & -10.751 & -10.439 & 5TF6 & -11.639 & -10.553 \\
5SZE & -10.459 & -9.448 & 5WZG & -10.218 & -10.425 & 5M8I & -5.736 & -8.523 \\
5UDZ & -9.592 & -9.424 & 5YTV & -8.043 & -7.575 & 6FQ3 & -8.015 & -9.025 \\
5YTX & -8.007 & -7.615 & 5YKI & -12.030 & -10.131 & 6DCL & -10.647 & -9.099 \\
6GD3 & -5.191 & -7.380 & 6G2K & -6.427 & -6.590 & 6CMN & -11.727 & -9.537 \\
5YTS & -7.577 & -7.593 & 5WWG & -9.077 & -8.937 & 5YTT & -7.394 & -7.684 \\
6FQR & -6.073 & -7.266 & 6GX6 & -7.264 & -6.113 & 5WWE & -8.550 & -8.960 \\
6GC5 & -8.435 & -6.263 & 6RT6 & -6.755 & -7.768 & 6NOF & -10.617 & -10.414 \\
6R7B & -10.558 & -9.604 & 6NOC & -10.490 & -10.480 & 6NOH & -10.421 & -10.600 \\
6A6J & -9.068 & -8.690 & 6NOD & -10.280 & -11.608 & 6G99 & -6.176 & -8.814 \\
6RT7 & -7.405 & -7.768 & 6U9X & -10.086 & -9.985 & 6NY5 & -9.952 & -10.305 \\
6GBM & -5.556 & -9.940 & 1EC6 & -7.990 & -9.086 & 1M8Y & -9.270 & -12.676 \\
1RPU & -13.300 & -9.964 & 1UTD & -16.890 & -9.426 & 2B6G & -10.600 & -9.984 \\
2ERR & -12.200 & -8.837 & 2F8K & -10.300 & -9.838 & 2G4B & -7.830 & -8.527 \\
2KFY & -8.840 & -8.242 & 2KG0 & -7.400 & -8.542 & 2KG1 & -7.450 & -8.217 \\
2KX5 & -11.400 & -9.900 & 2L41 & -4.250 & -7.366 & 2LA5 & -11.500 & -9.941 \\
2LEB & -9.320 & -9.401 & 2LEC & -9.440 & -9.318 & 2M8D & -8.730 & -8.704 \\
2MJH & -9.670 & -8.855 & 2MXY & -7.900 & -6.876 & 2MZ1 & -6.870 & -7.895 \\
2N82 & -10.300 & -10.288 & 2RU3 & -8.520 & -9.046 & 2XC7 & -7.040 & -8.838 \\
2XFM & -8.240 & -7.622 & 2XNR & -5.400 & -7.598 & 2ZKO & -8.050 & -10.387 \\
3BSB & -10.200 & -9.906 & 3BSX & -11.500 & -10.581 & 3BX2 & -9.970 & -9.808 \\
3K49 & -12.050 & -10.919 & 3K4E & -10.800 & -11.756 & 3K5Q & -9.300 & -8.994 \\
3K5Y & -10.540 & -9.156 & 3K5Z & -9.380 & -9.929 & 3K61 & -8.890 & -9.441 \\
3K62 & -8.740 & -9.403 & 3K64 & -9.220 & -9.071 & 3L25 & -8.175 & -9.369 \\
3LQX & -12.220 & -10.865 & 3MDG & -7.850 & -8.561 & 3MOJ & -13.000 & -10.732 \\
3NCU & -10.300 & -9.208 & 3NMR & -8.440 & -7.498 & 3NNH & -6.190 & -8.508 \\
3O3I & -6.520 & -7.067 & 3O6E & -7.070 & -6.518 & 3Q0L & -12.300 & -11.788 \\
3Q0M & -12.500 & -12.841 & 3Q0N & -10.500 & -12.090 & 3Q0P & -12.800 & -10.631 \\
3Q0Q & -13.300 & -12.998 & 3Q0R & -13.800 & -12.026 & 3Q0S & -10.800 & -13.101 \\
3QG9 & -10.600 & -11.421 & 3U4M & -15.900 & -10.450 & 3V71 & -11.100 & -10.349 \\
3V74 & -11.600 & -11.450 & 3WBM & -9.980 & -10.604 & 4CIO & -9.760 & -8.548 \\
4ED5 & -9.080 & -8.982 & 4ERD & -9.540 & -9.967 & 4JVH & -9.590 & -8.894 \\
4KJI & -8.940 & -10.456 & 4LG2 & -9.540 & -8.543 & 4NL3 & -7.380 & -10.121 \\
4O26 & -8.200 & -11.407 & 4OE1 & -12.110 & -9.565 & 4QI2 & -9.440 & -10.179 \\
4QVC & -11.560 & -10.097 & 4QVD & -10.120 & -10.982 & 4R3I & -7.770 & -6.926 \\
4RCJ & -8.180 & -8.311 & 4TUX & -10.400 & -9.742 & 4U8T & -9.130 & -7.530 \\
4Z31 & -8.790 & -8.985 & 5DNO & -7.830 & -8.605 & 5ELR & -6.250 & -8.884 \\
5V7C & -6.990 & -8.962 & 5W1I & -12.790 & -11.236 & 5WZH & -10.460 & -10.608 \\
6D12 & -9.300 & -10.239 & & &  & & &  \\
	
	\bottomrule
\end{longtable}
\newpage		
	